\documentclass{article}

\usepackage{microtype}
\usepackage{graphicx}
\usepackage{subfigure}
\usepackage{booktabs} 
\usepackage{multicol}

\usepackage{hyperref}

\usepackage[accepted]{icml2020}

\icmltitlerunning{Training Linear Neural Networks}

\usepackage{amsthm}
\usepackage{amsmath}
\usepackage{amssymb}
\usepackage{enumitem}

\def\R{\mathbb{R}}

\def\l{\left}
\def\r{\right}
\def\a{\alpha}
\def\d{\delta}

\def\der{\text{d}}
\def\b{\beta}

\def\rank{\text{rank}}
\def\A{\mathcal{A}}
\def\M{\mathcal{M}}
\def\ol{\overline}

\def\range{\text{range}}
\def\s{\sigma}
\def\D{\Delta}

\def\g{\gamma}
\def\neigh{\mathcal{N}}
\def\sign{\text{sign}}
\def\svd{\text{tSVD}}
\def\wt{\widetilde}
\def\rowspan{\text{row span}}
\def\proj{\mathcal{P}}

\def\NB{\mathbf{N}}
\def\th{\text{th}}
\def\null{\text{null}}

\def\range{\text{range}}
\def\map{\Pi}

\newtheorem{thm}{Theorem}[section]
\newtheorem{lem}[thm]{Lemma}

\newtheorem{prop}[thm]{Proposition}
\newtheorem{assumption}[thm]{Assumption}

\newtheorem{defn}[thm]{Definiton}

\usepackage{tikz}
\newcommand*\circled[1]{\tikz[baseline=(char.base)]{
             \node[shape=circle,draw,inner sep=.8pt] (char) {#1};}}

\begin{document}

\twocolumn[
\icmltitle{Training Linear Neural Networks:\\ 
Non-Local Convergence and Complexity Results
}




\begin{icmlauthorlist}
\icmlauthor{Armin Eftekhari}{to}
\end{icmlauthorlist}

\icmlaffiliation{to}{Department of Mathematics and Mathematical Statistics, Umea University, Sweden. AE is indebted to Holger Rauhut, Ulrich Terstiege and Gongguo Tang for insightful discussions}
\icmlcorrespondingauthor{Armin Eftekhari}{armin.eftekhari@umu.se}

\icmlkeywords{Linear networks, lazy training, non-strict saddle points, gradient flow}

\vskip 0.3in
]



\printAffiliationsAndNotice{}  

\begin{abstract}
Linear networks provide  valuable insights into the workings of neural networks in general. 
This paper identifies conditions under which the gradient flow provably trains a linear network, in spite of the non-strict saddle points present in the optimization landscape. 
This paper also provides the computational complexity of training linear networks with gradient flow.
To achieve these results, this work develops a  machinery to provably identify the stable set of gradient flow, which then enables us to
improve over the state of the art in the literature of linear networks~\cite{bah2019learning,arora2018convergence}.
Crucially, our results appear to be the first to break away from the {lazy training} regime which has dominated the literature of neural networks. 
This work requires the network to have \emph{a} layer with one neuron, which subsumes the  networks  with a scalar output, but extending the results of this theoretical work to all linear networks remains a challenging open problem.  \end{abstract}

\section{Introduction and Overview}\label{sec:intro}

Consider the training samples and their labels 
$
\{x_i,y_i\}_{i=1}^m \subset \R^{d_x}\times  \R^{d_y}$, respectively. 
By concatenating $\{x_i\}_i$ and $\{y_i\}_i$, we form the data  matrices
\begin{align}
X\in \R^{d_x\times m}, \qquad Y\in \R^{d_y\times m}.
\label{eq:trainigData}
\end{align}
Consider also a linear network, i.e., 
a neural network where the nonlinear activation functions are replaced with the identity map. 
To be specific, 
with $N$ layers and the corresponding weight matrices $\{W_i\}_{i=1}^N$, 
this network is characterized by the linear map
\begin{align}
    & \R^{d_x} \rightarrow \R^{d_y}\nonumber\\
    & x\rightarrow Wx,
    \label{eq:linNet}
\end{align}
and the matrix $W\in \R^{d_y\times d_x}$ in~\eqref{eq:linNet} is itself specified with the nonlinear (and often over-parametrized) map 
\begin{align}
     & \R^{d_1\times d_0} \times \cdots \times \R^{d_N \times d_{N-1}} \longrightarrow \R^{d_y\times d_x}  \nonumber\\
    & (W_1,\cdots,W_N) \longrightarrow W:= W_N \cdots W_1,
    \label{eq:linearNetMap}
\end{align}
 where we set $d_0 = d_x$ and  $d_N = d_y$  for consistency.

In foregoing the full generality of nonlinear neural networks, linear networks afford us a level of insight and technical rigor that is out of the reach for  nonlinear networks, at least with our current technical tools~\cite{arora2018optimization,yan1994global,kawaguchi2016deep,chitour2018geometric,trager2019pure,saxe2013exact,lu2017depth,yun2017global}.

Indeed, despite the absence of activation functions,  matrix $W$ in the linear network~(\ref{eq:linNet},\ref{eq:linearNetMap}) is a nonlinear function of $\{W_i\}_i$, and  training this linear network thus involves solving a  nonconvex optimization problem in $\{W_i\}_i$, which shares many interesting features of the nonlinear neural networks.
Simply put,  we cannot claim to understand neural networks in general without  understanding linear  
networks.

Training the linear network~(\ref{eq:linNet},\ref{eq:linearNetMap}) with the  data $(X,Y)$ in~\eqref{eq:trainigData} can be cast as the optimization problem 
\begin{align}
\begin{cases}    
    \underset{{W_1,\cdots,W_N}}{\min}\,\, \frac{1}{2}\| Y - W_N W_{N-1} \cdots W_1 X\|_F^2 \\
     \text{subject to } W_j\in \R^{d_j\times d_{j-1}} \qquad \forall j\in [N],
\end{cases}
\label{eq:mainN0}
\end{align}
which is nonconvex when $N\ge 2$. Above, $[N]=\{1,\cdots,N\}$. 
Let us introduce  the shorthand 
\begin{align}
 W_{\NB}& :=(W_1,\cdots,W_N)  \nonumber\\
& \in \R^{d_1\times d_0} \times \cdots \times \R^{d_N\times d_{N-1}}=:\R^{d_{\NB}},
\label{eq:sharthandWNWj}
\end{align}
which allows us to rewrite problem~\eqref{eq:mainN0} more compactly as 
\begin{align}
    \underset{W_{\NB}}{\min}\, \, L_N(W_{\NB})\,\,\text{subject to}\,\, W_{\NB}\in \R^{d_{\NB}},
\label{eq:mainN}
\end{align}
where $L_N(W_{\NB}) := \frac{1}{2}\|Y-W_N\cdots W_1X\|_F^2$. 
With this setup and before turning to the details, let us highlight the contributions of this paper, in the order of appearance.
\begin{itemize}[leftmargin=*]
    \item Theorem~\ref{thm:wellBehavedMain} in Section~\ref{sec:landscape} provides a new analysis of the  optimization landscape of linear networks, where we uncover a previously-unknown link to the celebrated Eckart-Young-Mirsky theorem and the geometry of the principal component analysis (PCA).
    
    \item Theorem~\ref{cor:conjProved} in Section~\ref{sec:cvg}  identifies the conditions under which gradient flow successfully trains a linear network, despite the presence of non-strict saddle points in the optimization landscape.
    
    Theorem~\ref{cor:conjProved} thus improves  the state of the art in~\cite{bah2019learning} as the first convergence result outside of the {lazy training} regime, reviewed later.  
    This improvement is achieved with a new argument that provably identifies the stable set of gradient flow, in the language of dynamical systems theory.
    
    Theorem~\ref{cor:conjProved} applies to linear networks that have \emph{a} layer with a single neuron, see Assumption~\ref{assump:key2}. This case corresponds to the popular {spiked covariance} model in statistics, and subsumes networks with a scalar output.  Extension of Theorem~\ref{cor:conjProved} to all linear networks remains a challenging open problem, as there is no natural notion of stable set in general. We will however speculate about how Theorem~\ref{cor:conjProved} might serve as the natural building block for a more general  result in the future.

    \item  Theorem~\ref{thm:MainSimple} in Section~\ref{sec:networkDepth} quantifies the computational complexity of training a linear network by   establishing \emph{non-local} convergence rates for gradient flow. This result  also quantifies for the first time how the (faraway) convergence rate benefits from increasing the network depth.
    
    Theorem~\ref{thm:MainSimple} presents the first convergence rate for linear networks outside of the {lazy training} regime, thus improving the state of the art in~\cite{arora2018convergence}. 
    
    Indeed, of the dozens of related works, virtually all belong to this lazy regime, to the best of our knowledge, thus signifying the importance of this breakthrough.
    Theorem~\ref{thm:MainSimple} also corresponds to the spiked covariance model. 
\end{itemize}

\section{Landscape of Linear Networks}\label{sec:landscape}

The landscape of  nonconvex program~\eqref{eq:mainN} has been widely studied in the literature, with  contributions from~\cite{arora2018optimization,arora2018convergence,chitour2018geometric,bartlett2019gradient,saxe2013exact,hardt2016identity,laurent2018deep,trager2019pure,baldi1989neural,zhu2019global,he2016identity,nguyen2019connected}. 
The state of the art here is Proposition~31 in~\cite{bah2019learning}, reviewed in Appendix~\ref{sec:anotherProofBah}, which itself  improves over Theorem~2.3 in~\cite{kawaguchi2016deep}.

The main result of this section, Theorem~\ref{thm:wellBehavedMain} below, is a variation of Proposition~31 in~\cite{bah2019learning} with an additional assumption. In Section~\ref{sec:cvg}, we will use Theorem~\ref{thm:wellBehavedMain} to improve  the state of art for training linear networks, under this new assumption. 

Crucially, the proof of Theorem~\ref{thm:wellBehavedMain} uncovers a  previously-unknown link to the celebrated Eckart-Young-Mirsky theorem~\cite{eckart,mirsky} and the geometry of PCA.






To begin, let us concretely define the notion of optimality for problem~\eqref{eq:mainN}.   
\begin{defn}[First-order stationarity for (\ref{eq:mainN})] \label{defn:fosp}
We say that $\ol{W}_{\NB}\in \R^{d_{\NB}}$ is a first-order stationary point (FOSP) of problem~\eqref{eq:mainN} if 
\begin{align}
    \nabla L_N(\ol{W}_{\NB}) = 0, \label{eq:firstOrderLN}
\end{align}
where $\nabla L_N(\ol{W}_{\NB})$ is the gradient of $L_N$ at $\ol{W}_{\NB}$. 
\end{defn}

\begin{defn}[Second-order stationarity for (\ref{eq:mainN})]\label{defn:secondStN}
We say that $\ol{W}_{\NB}\in \R^{d_{\NB}}$ is a second-order stationary point (SOSP) of problem~\eqref{eq:mainN} if, in addition to~\eqref{eq:firstOrderLN}, it  holds that 
\begin{align}
    \nabla^2 L_N(\ol{W}_{\NB})[\D_{\NB}] \ge  0,\qquad 
    \forall \D_{\NB}\in \R^{d_{\NB}},
    \label{eq:secondOrderLN}
\end{align}
where  $\nabla^2L_N (\ol{W}_{\NB})[\D_{\NB}]$ is the  second derivative of $L_N$  at $W_{\NB}$ along the direction $\D_{\NB}$.
\end{defn}
\begin{defn}[Strict saddles of (\ref{eq:mainN})]\label{defn:strictSaddle}
Any FOSP of problem~\eqref{eq:mainN}, which is not an SOSP, is a strict saddle point of~\eqref{eq:mainN}.
\end{defn}

Any SOSP of problem~\eqref{eq:mainN} is either a local minimizer of~\eqref{eq:mainN}, or a non-strict saddle point  of problem~\eqref{eq:mainN}. Unlike a non-strict saddle point, there always exists a descent direction to escape from a strict saddle point~\cite{lee2017first}.
To continue, let 
\begin{align}
r := \min_{j\le N} d_j,\qquad \text{(smallest width of the network)}
\label{eq:rankNet}
\end{align}
denote the smallest width of the linear network~(\ref{eq:linNet},\ref{eq:linearNetMap}). 
As shown in Appendix~\ref{sec:surjective}, we can reformulate problem~\eqref{eq:mainN} as
\begin{subequations}
\begin{align} 
     & 
     \underset{W_{\NB}}{\min} \,\, L_N(W_{\NB})
     \nonumber\\
     & = 
    \begin{cases}
     \underset{W}{\min} \,\, \frac{1}{2}\| Y - W X\|_F^2=: L_1(W)\\
     \text{subject to } \rank(W) \le r
    \end{cases} \label{eq:mainOne} \\
    & = \begin{cases}
     \underset{P,Q}{\min} \,\, \frac{1}{2}\| Y - PQ X\|_F^2=: L_2(P,Q) \\
     \text{subject to } P\in \R^ {d_y\times r},\, Q\in \R^{r\times d_x}. \label{eq:mainTwo}
    \end{cases}
\end{align}
\end{subequations}
In particular, the notion of optimality for problem~\eqref{eq:mainTwo}  is defined similar to Definitions~\ref{defn:fosp}-\ref{defn:strictSaddle}.
There is a correspondence between the stationary points of problems~\eqref{eq:mainN} and~\eqref{eq:mainTwo},   proved in Appendix~\ref{sec:proofLNtoL2}.

\begin{lem}[Pairwise correspondence between SOSPs]\label{lem:LNtoL2}
Any FOSP $\ol{W}_{\NB}=(\ol{W}_1,\cdots,\ol{W}_N)$ of problem~\eqref{eq:mainN}  corresponds to an FOSP $(\ol{P},\ol{Q})$ of problem~\eqref{eq:mainTwo}, provided that $\ol{W} = \ol{W}_N\cdots \ol{W}_1$ is rank-$r$. 
Moreover, any SOSP $\ol{W}_{\NB}$ of problem~\eqref{eq:mainN} corresponds to an SOSP $(\ol{P},\ol{Q})$ of problem~\eqref{eq:mainTwo}, 
provided that $\rank(\ol{W}) = r$.

\end{lem}

Let $\proj_X:=X^\dagger X$ and $\proj_{X^\perp}:=I_m - \proj_X$ denote the orthogonal projections onto the row span of $X$ and its orthogonal complement, respectively. Here, $X^\dagger$ is the pseudo-inverse of $X$ and $I_m\in \R^{m\times m}$ is the identity matrix. 
Using the decomposition $Y=Y\proj_X+Y\proj_{X^\perp}$, we can in turn rewrite  problem~\eqref{eq:mainTwo} as 
\begin{align}
& \min_{P,Q} \,\, L_2(P,Q) \nonumber\\
  & = \frac{1}{2}\|Y\proj_{X^\perp}\|_F^2
  + \begin{cases}
    \underset{P,Q,Q'}{\min} \,\, \frac{1}{2}\|Y\proj_X-PQ'\|_F^2 \nonumber\\
    \text{subject to } Q' = Q X 
    \end{cases} \nonumber\\
    & \ge \frac{1}{2}\|Y\proj_{X^\perp}\|_F^2
  + 
    \underset{P,Q'}{\min}\,\, \frac{1}{2}\|Y\proj_X-PQ'\|_F^2.
\end{align}
The relaxation above is tight, and  there is a correspondence between the stationary points, proved in Appendix~\ref{sec:proofL2toPCA}.
\begin{lem}[Pairwise correspondence between SOSPs]\label{lem:L2toPCA}
Suppose that $XX^\top$ is invertible. Then it holds that
\begin{align}
    & -\frac{1}{2}\|Y\proj_{X^\perp}\|_F^2+
    \underset{P,Q}{\min}\,\, L_2(P,Q) \nonumber\\
    & =  
    \begin{cases}
     \underset{P,Q'}{\min} \,\, \frac{1}{2}\|Y\proj_X - PQ'\|_F^2=: L_2'(P,Q') \\
    \emph{\text{subject to }} P\in \R^{d_y\times r}, \, Q'\in \R^{r\times m}.
    \end{cases}
        \label{eq:pca}
\end{align}
Any FOSP $(\ol{P},\ol{Q})$ of problem~\eqref{eq:mainTwo} corresponds to an FOSP~$(\ol{P},\ol{Q}')$ of problem~\eqref{eq:pca}.
Moreover, any SOSP $(\ol{P},\ol{Q})$ of problem~\eqref{eq:mainTwo} corresponds to an SOSP $(\ol{P},\ol{Q}')$ of problem~\eqref{eq:pca}. 
\end{lem}

Note that solving problem~\eqref{eq:pca} involves finding a best rank-$r$ approximation of $Y\proj_X$ or, equivalently, finding $r$ leading principal components of $Y\proj_X$~\cite{murphy2012machine}.   
By combining Lemmas~\ref{lem:LNtoL2} and~\ref{lem:L2toPCA}, we immediately reach the following conclusion.
\begin{lem}[Pairwise correspondence between SOSPs]\label{lem:cmobinedLand}
Suppose that $XX^\top$ is invertible. Then  any FOSP $\ol{W}_{\NB}$ of problem~\eqref{eq:mainN} corresponds to an FOSP $(\ol{P},\ol{Q}')$ of problem~\eqref{eq:pca}, provided that $\ol{W}=\ol{W}_N\cdots \ol{W}_1$ is rank-$r$. 
Moreover, any SOSP $\ol{W}_{\NB}$ of problem~\eqref{eq:mainN} corresponds to an SOSP~$(\ol{P},\ol{Q}')$ of~\eqref{eq:pca}, 
provided that $\rank(\ol{W})=r$.
\end{lem}

We next recall a variant of the celebrated EYM theorem~\cite{eckart,mirsky,hauser2018pca,hauser2018pca2}, which specifies the landscape of the PCA problem~\eqref{eq:pca}. 
\begin{thm}[EYM theorem]\label{thm:eym} Any SOSP $(\ol{P},\ol{Q}')$ of the PCA problem~\eqref{eq:pca} is also a global minimizer of problem~\eqref{eq:pca}. 
\end{thm}

With Lemma~\ref{lem:cmobinedLand} at hand, we invoke Theorem~\ref{thm:eym}  to uncover the landscape of problem~\eqref{eq:mainN}, see Appendix~\ref{sec:proofThmMain} for the proof.

\begin{thm}[Landscape of linear networks]\label{thm:wellBehavedMain}
    Suppose that $XX^\top$ is invertible and that $\rank(Y X^\dagger X)\ge r$. 
    Then any SOSP $\ol{W}_{\NB}=(\ol{W}_1,\cdots,\ol{W}_N)$ of problem~\eqref{eq:mainN}  is a global minimizer of problem~\eqref{eq:mainN}, provided that $\ol{W}_N\cdots \ol{W}_1$ is a rank-$r$ matrix.  
\end{thm}

In words, Theorem~\ref{thm:wellBehavedMain} identifies certain SOSPs of problem~\eqref{eq:mainN} which are  global minimizers of problem~\eqref{eq:mainN}. 
A few important remarks are in order: 
{\circled{1}}~The proof of Theorem~\ref{thm:wellBehavedMain} establishes a pairwise correspondence with the stationary points and the geometry of the PCA problem. This connection was previously unknown, to the best of our knowledge. 

{\circled{2}}~Any rank-degenerate SOSP of problem~\eqref{eq:mainN} is  excluded from Theorem~\ref{thm:wellBehavedMain}, i.e., any SOSP $\ol{W}_{\NB}=(\ol{W}_1,\cdots,\ol{W}_N)$ such that $\rank(\ol{W}_N\cdots \ol{W}_1)<r$ is excluded from our result. 
For example, the zero matrix is a {spurious} SOSP 
of problem~\eqref{eq:mainN} when the network depth $N\ge 2$, as observed in~\cite{bah2019learning,kawaguchi2016deep,trager2019pure}. 
 The landscape of problem~\eqref{eq:mainN} in general is therefore more complicated than the special case of $N=2$, corresponding to the Eckart-Young-Mirsky theorem, see Theorem~\ref{thm:eym}. 
 
{\circled{3}}~
Theorem~\ref{thm:wellBehavedMain} is a variation of Proposition~31 in~\cite{bah2019learning} with a new assumption on $Y X^\dagger X$, which  will be necessary shortly.
Similar assumptions have been used in the context of PCA, for example in~\cite{helmke1995critical}.  
{\circled{4}}~For completeness, we also prove Theorem~\ref{thm:wellBehavedMain} using Proposition~31 in~\cite{bah2019learning} as the starting point, see Appendix~\ref{sec:anotherProofBah}. 

\section{Convergence of Gradient Flow}\label{sec:cvg}

In view of Theorem~\ref{thm:wellBehavedMain} above, even though nonconvex, the landscape of problem~\eqref{eq:mainN}  has certain favourable properties. On the other hand, problem~\eqref{eq:mainN}  fails to satisfy the {strict saddle} property that enables first-order algorithms to avoid saddle points~\cite{lee2016gradient,ge2015escaping}. For example, the zero matrix is a non-strict saddle point of problem~\eqref{eq:mainN} when the network depth $N\ge 2$, as discussed earlier.

Against this mixed background, it is natural to ask if first-order  algorithms can successfully train linear neural networks. This fundamental question has remained unanswered in the  literature, to our knowledge.
Indeed, outside the {lazy training} regime, 
reviewed in Section~\ref{sec:networkDepth}, it is not known if gradient flow can successfully solve problem~\eqref{eq:mainN}.


In fact, the state of the art here, Theorem~35(a) in~\cite{bah2019learning},  guarantees the convergence of gradient flow to a  minimizer of $L_N$, when \emph{restricted} to one of few regions in $\R^{d_{\NB}}$. Even though these regions are known in advance, their result cannot predict which region  would contain the limit point of gradient flow, for a given initialization.  

In other words, Theorem~35(a) in~\cite{bah2019learning} does  \emph{not} guarantee the convergence of gradient flow to a global minimzer of problem~\eqref{eq:mainN}, and   gradient flow might indeed converge to a spurious SOSP of problem~\eqref{eq:mainN}, such as the zero matrix. {For completeness, Theorem~35(a) in~\cite{bah2019learning} is reviewed in
Appendix~\ref{sec:revBahMain}.
}

 In an important setting, this section answers the open question of convergence of gradient flow for training linear networks. {This is achieved with a new argument, which enables us to provably identify the {stable set} of gradient flow, in the language of dynamical systems theory.}
 Let us begin with the necessary preparations.

Consider gradient flow applied to program \eqref{eq:mainN}, specified as
\begin{align}
     & \dot{W}_j(t)= \frac{\der W_j(t)}{\der t} = - \nabla_{W_j} L_N\l(W_{\NB}(t)\r), \nonumber\\ 
    & \forall j\in [N], \, \forall t\ge 0,
    \qquad \text{(gradient flow)}
    \label{eq:gradFlowN}
\end{align}
and initialized at $W_{\NB,0}\in \R^{d_{\NB}}$. Above, $\nabla_{W_j} L_N$ is the gradient of $L_N$ with respect to  $W_j$, the weight matrix for the $j^{\th}$ layer of the linear network, see~(\ref{eq:mainN0},\ref{eq:sharthandWNWj},\ref{eq:mainN}).

A consequence of the Lojasiewicz' theorem
is the following convergence result for  the gradient flow.  See  Appendix~\ref{sec:proofGfCvgFOSP} for the proof, similar to Theorem~11 in~\cite{bah2019learning}. 

\begin{lem}[Convergence, uninformative]\label{lem:gfCvgFOSP}
If $XX^\top$ is invertible, then gradient flow~\eqref{eq:gradFlowN} converges.
Moreover, the limit point is an SOSP $\ol{W}_{\NB}\in \R^{d_{\NB}}$ of problem~\eqref{eq:mainN}, for almost every initialization $W_{\NB,0}$ with respect to the Lebesgue measure in $\R^{d_{\NB}}$.  
\end{lem}

To study  this limit point $\ol{W}_{\NB}$, we focus here on a common initialization technique for linear networks~\cite{hardt2016identity,bartlett2019gradient,arora2018convergence,arora2018optimization}.

\begin{defn}[Balanced initialization]\label{defn:balanced}
For gradient flow~\eqref{eq:gradFlowN}, we call $W_{\NB,0}=(W_{1,0},\cdots,W_{N,0})\in \R^{d_{\NB}}$ a balanced initialization if 
\begin{align}
    & W_{j+1,0}^\top W_{j+1,0} = W_{j,0} W_{j,0}^\top,
    \qquad \forall j\in [N-1].\label{eq:balanced}
\end{align}
\end{defn} 
Claim~4 in~\cite{arora2018convergence} underscores the necessity of a (nearly) balanced initialization for linear networks. More generally,~\cite{sutskever2013importance} highlights the importance of initialization in deep neural networks.  
The main result of this section, 
Theorem~\ref{cor:conjProved} below, thus requires an (exactly) balanced initialization.  

Assuming exact balanced-ness in~\eqref{eq:balanced} is sufficient here because the focus of this theoretical work is continuous-time optimization. More generally, approximate balanced-ness is necessary for discretized algorithms, such as gradient descent~\cite{arora2018convergence}. We avoid this additional layer of complexity here as it does not seem to add any key theoretical  insights to this paper.

 A useful observation is that, 
if the initialization is balanced, gradient flow  remains balanced afterwards, see for example Lemma~2 in~\cite{bah2019learning}. 
More formally, gradient flow~\eqref{eq:gradFlowN} satisfies
\begin{align}
    & W_{j+1,0}^\top W_{j+1,0} = W_{j,0} W_{j,0} ^\top, \qquad \forall j\in [N-1] 
  \nonumber
    \\
    & \Longrightarrow W_{j+1}(t)^\top W_{j+1}(t) = W_{j}(t) W_{j}(t)^\top, \label{eq:remainBalanced}
\end{align}
for every $j\in [N-1]$ and every $t\ge 0$. Above, the weight matrix $W_j(t)$ is the $j^{\th}$ component of $W_{\NB}(t)$, see~\eqref{eq:sharthandWNWj}.

Alongside gradient flow~\eqref{eq:gradFlowN}, it is convenient to introduce another flow~\cite{bah2019learning,arora2018convergence}, which dictates the evolution of the end-to-end product of the weight matrices of the linear network.

Concretely, for a 
matrix $W\in \R^{d_y\times d_x}$,  consider the linear operator $\A_W$ specified as 
\begin{align}
    \A_W :& \R^{d_y\times d_x} \rightarrow \R^{d_y\times d_x} \nonumber\\
    & \D \rightarrow \sum_{j=1}^N (W W^\top)^{\frac{N-j}{N}} \D (W^\top W)^{\frac{j-1}{N}}.
    \label{eq:defnAW}
\end{align}
For a  balanced initialization $W_{\NB,0}=(W_{1,0},\cdots,W_{N,0})$, gradient flow~\eqref{eq:gradFlowN} in $\R^{d_{\NB}}$  induces a flow in  $\R^{d_y\times d_x}$, initialized at $W_0 = W_{N,0}\cdots W_{1,0}\in\R^{d_y\times d_x}$ and  specified as  
\begin{align}
     \dot{W}(t)&  = - \A_{W(t)}\l( \nabla L_1 (W(t)) \r) \qquad \forall t\ge 0, \nonumber\\
     & = - \A_{W(t)}( W(t) XX^\top - YX^\top ) 
     \qquad \text{(see \eqref{eq:mainOne})}
     \nonumber\\
    & \qquad \qquad   \text{(induced flow)}
    \label{eq:flowW}
\end{align}
see for example Equation~(26) in~\cite{bah2019learning}. Above, 
\begin{align}
W(t)  = W_N(t)\cdots W_1(t)\in \R^{d_y\times d_x}.
\label{eq:gradFlowToIndFlow}
\end{align}
We will refer to~\eqref{eq:flowW}  as the \emph{induced flow}, which governs the evolution of the end-to-end product of the weight matrices of the linear network.

It is  known that  induced flow~\eqref{eq:flowW} admits an analytic singular value decomposition (SVD), see for example Lemma~1 and Theorem~3 in~\cite{arora2019implicit} or~\cite{analyticLecturess}. More specifically, it holds that 
\begin{align}
    W(t) \overset{\text{SVD}}{=} \wt{U}(t) \wt{S}(t) \wt{V}(t)^\top,\qquad  \forall t\ge 0, 
    \label{eq:analyticSVDCVG}
\end{align}
provided that the network depth $N\ge2$. 
In~\eqref{eq:analyticSVDCVG}, $\wt{U}(t),\wt{V}(t),\wt{S}(t)$ are analytic functions of $t$~\cite{parks1992primer}. Moreover,
$\wt{U}(t)\in \R^{d_y\times d_y}, \wt{V}(t)\in \R^{d_x\times d_x}$ are orthonormal bases, and $\wt{S}(t)\in \R^{d_y\times d_x}$ contains the singular values of $W(t)$ in no specific order.

The evolution of the singular values of $W(t)$ in~\eqref{eq:analyticSVDCVG} is also known~\cite{arora2019implicit,townsend2016differentiating}. In particular, the following byproduct about the rank of $W(t)$ is important for us, see Appendix~\ref{sec:proofRankCteFlow} for the proof. 
\begin{lem}[Rank-invariance]\label{lem:rankCteFlow}
For induced flow~\eqref{eq:flowW}, $\rank(W(t))=\rank(W_0)$ for all $t\ge 0$, provided that  $XX^\top$ is invertible and  the network depth $N\ge 2$. 
\end{lem}

Let us henceforth assume that $XX^\top$ is invertible, and that gradient flow~\eqref{eq:gradFlowN} is initialized at $W_{\NB,0}\in \M_{\NB,r}$, where 
\begin{align}
     \M_{\NB,r} & := 
    \Big\{ 
    W_{\NB} : \rank(W_N \cdots W_1) = r
    \Big\}
     \subset \R^{d_{\NB}},
    \label{eq:preimage}
\end{align}
see~\eqref{eq:sharthandWNWj}.
We now make the following observations about the set $\M_{\NB,r}$, proved in Appendix~\ref{sec:proofGenericNet}.

\begin{lem}[Propeties of $\M_{\NB,r}$]\label{lem:genericNet} \emph{\circled{1}} $\M_{\NB,r}$ is not a closed subset of $\R^{d_{\NB}}$. \emph{\circled{2}}~The complement of $\M_{\NB,r}$ in $\R^{d_{\NB}}$ has Lebesgue measure zero. (In particular, $\M_{\NB,r}$ is a dense subset of $\R^{d_{\NB}}$.) 
\end{lem}

In view of Lemma~\ref{lem:genericNet}, almost every initialization $W_{\NB,0}\in\R^{d_{\NB}}$ of gradient flow~\eqref{eq:gradFlowN} falls into the set $\M_{\NB,r}$, i.e.,
\begin{align}
    W_{\NB,0}\in\M_{\NB,r}, \qquad \text{almost surely}.
    \label{eq:asInM}
\end{align}
Moreover, once initialized  in $\M_{\NB,r}$ with a balanced initialization, induced flow~\eqref{eq:flowW} remains rank-$r$ at all times by~(\ref{eq:gradFlowToIndFlow},\ref{eq:preimage}) and Lemma~\ref{lem:rankCteFlow}. Consequently, gradient flow~\eqref{eq:gradFlowN} remains in $\M_{\NB,r}$ at all times, see again~(\ref{eq:gradFlowToIndFlow},\ref{eq:preimage}). We combine this last observation with~\eqref{eq:asInM} to conclude that
\begin{align}
    W_{\NB}(t) \in \M_{\NB,r}, \qquad \forall t\ge 0, \quad \text{ almost surely},
    \label{eq:remainsManifoldRdN}
\end{align}
over the choice of balanced initialization $W_{\NB,0}\in \R^{d_{\NB}}$. 
Despite~\eqref{eq:remainsManifoldRdN}, the limit point $\ol{W}_{\NB}$ of gradient flow~\eqref{eq:gradFlowN} might \emph{not} belong to $\M_{\NB,r}$ because $\M_{\NB,r}$ is \emph{not} closed, see Lemma~\ref{lem:genericNet}. 

That is, even though  the limit point  $\ol{W}_{\NB}$ of gradient flow is almost surely an SOSP of problem~\eqref{eq:mainN} by Lemma~\ref{lem:gfCvgFOSP}, we \emph{cannot} apply Theorem~\ref{thm:wellBehavedMain} and
 $\ol{W}_{\NB}$ might be an spurious SOSP of problem~\eqref{eq:mainN}, such as the zero matrix. 
Indeed, Remark~39 in~\cite{bah2019learning} constructs an example where $\ol{W}_{\NB}\notin \M_{\NB,r}$, see also~\cite{yan1994global}. To avoid this unwanted behaviour, it is  necessary to restrict the initialization of the gradient flow and impose additional assumptions.

Our first assumption is that the data is statistically whitened, which is common in the analysis of  linear networks~\cite{arora2018convergence,bartlett2019gradient}. 

\begin{defn}[Whitened data]\label{defn:whitened} We say that the data matrix $X\in \R^{d_x\times m}$ is whitened if 
\begin{align}
    \frac{XX^\top}{m} = \frac{1}{m}\sum_{i=1}^m x_i x_i^\top =  I_{d_x},
    \label{eq:whitenedData}
\end{align}
where $I_{d_x}\in \R^{d_x\times d_x}$ is the identity matrix. 
\end{defn}

Our second assumption is that $r=1$ in~\eqref{eq:rankNet}. This case is significant as it corresponds to the popular \emph{spiked covariance} model in statistics and signal  processing~\cite{eftekhari2019moses,johnstone2001distribution,vershynin2012close,berthet2013optimal,deshpande2014information}, to name  a few. 

Moreover, $r=1$ subsumes the important case of networks with a scalar output. 

Lastly, the case $r=1$  appears  to be the natural building block for the case $r>1$ via a \emph{deflation argument}~\cite{mackey2009deflation,zhang2011large}. Indeed,  gradient flow~\eqref{eq:gradFlowN} moves orthogonal to the principal directions that it has previously discovered or ``peeled''. Extending our results to the case $r>1$ remains a challenging open problem.

From~(\ref{eq:mainOne}) with $r=1$, recall that problem~\eqref{eq:mainN} for training a linear neural network is closely related to the problem
\begin{align}
    & 
    \underset{W}{\min} \,\, \frac{1}{2}\| Y \proj_X - WX \|_F^2\,\, \text{subject to}\,\, \rank(W) \le r= 1
    \nonumber\\
    & = 
    \underset{W}{\min} \,\,  \frac{m}{2}\| Z - W\|_F^2\,\,  \text{subject to } \rank(W)\le 1,
    \label{eq:whyPCA}
\end{align}
where the second line above is obtained using~\eqref{eq:whitenedData}, and
\begin{align}
    Z:= \frac{YX^\top}{m}.
    \label{eq:defnZ}
\end{align}
We are in position to  collect  all the assumptions made in this section in one place.
\begin{assumption}\label{assump:key2} In this section, we assume that the linear network~\emph{(\ref{eq:linNet},\ref{eq:linearNetMap})} has depth $N\ge 2$,  and one of the layers has only one neuron, i.e.,  $r=1$ in~\eqref{eq:rankNet}.
Moreover, the data matrix $X$ in~\eqref{eq:trainigData} is whitened as in~\eqref{eq:whitenedData}, and  $Z=\frac{1}{m} YX^\top$ in~\eqref{eq:defnZ} satisfies  
\begin{align}
\rank(Z)\ge r= 1, \quad \g_Z:= \frac{s_{Z,2}}{s_{Z}}<1, 
\label{eq:simpler}
\end{align}
where $s_Z$ and $s_{Z,2}$ are the two largest singular values of~$Z$.
Lastly, we assume that the initialization  of gradient flow~\eqref{eq:gradFlowN} is balanced, see Definition~\ref{defn:balanced}. 
\end{assumption}

In view of~(\ref{eq:simpler}), let us define  
\begin{align}
     Z_1 = u_Z \cdot s_Z\cdot v_Z^\top
     \qquad \text{(target matrix)}
    \label{eq:decomZ0}
\end{align} 
to be the best rank-$1$ approximation of $Z$, obtained via SVD. 
Here, $\|u_Z\|_2=\|v_Z\|_2=1$, and $s_Z$ appeared in~\eqref{eq:simpler}.

Note that $Z_1$ is the unique solution of problem~\eqref{eq:whyPCA}, because $Z$ has a nontrivial spectral gap in~\eqref{eq:simpler}, see 
for example Section 1 of~\cite{golub1987generalization}.

Let us fix $\a\in [\g_Z,1)$.  
To exclude the zero matrix as the limit point of gradient flow~\eqref{eq:gradFlowN},
the key is to restrict the initialization to a particular subset of the feasible set of problem~\eqref{eq:mainN} with $r=1$, specified as 
\begin{align}
   &  \neigh_{\NB,\a}  := \Big\{ W_{\NB}=(W_1,\cdots,W_N): \nonumber\\
   &  \qquad W_N \cdots W_1 \overset{\svd}{=} u_W \cdot s_W \cdot v_W^\top,   \nonumber\\
    &\,\,  s_W >  (\a-\g_Z) s_{Z}, \, u_W^\top Z_1 v_W> \a s_{Z}
    \Big\} \subset \R^{d_{\NB}},
    \label{eq:neighWNB}
\end{align}
where  $s_Z,\g_{Z}$ were defined in~\eqref{eq:simpler}. Above, $\svd$ stands for the thin SVD. A simple observation is that the set $\neigh_{\NB,\a}$ has infinite Lebesgue measure in $\R^{d_{\NB}}$.

Such restriction of the feasible set of problem~\eqref{eq:mainN}  is  \emph{necessary} as described earlier, see also  the negative example constructed in Remark~39 of~\cite{bah2019learning}.
Crucially, note that the end-to-end matrices in $\neigh_{\NB,\a}$ are positively correlated with $Z_1$, and also bounded away from the origin.

An important observation is that, once initialized in $\neigh_{\NB,\a}$, gradient flow~\eqref{eq:gradFlowN} avoids the zero matrix, see Appendix~\ref{sec:proofRemainNeighSimpleNBDomain}.

\begin{lem}[Stable set]\label{lem:remainNeighSimpleNBDomain}
For gradient flow~\eqref{eq:gradFlowN} initialized at $W_{\NB,0}\in \neigh_{\NB,\a}$, the limit point exists and 
satisfies $\ol{W}_{\NB}\in \M_{\NB,1}$.
Above, $\a\in (\g_Z,1)$, and  Assumption~\emph{\ref{assump:key2}} and its notation are in force, see also~\emph{(\ref{eq:preimage},\ref{eq:neighWNB})}.
\end{lem}

Combining Lemma~\ref{lem:remainNeighSimpleNBDomain} with Lemma~\ref{lem:gfCvgFOSP}, we  find that the limit point $\ol{W}_{\NB}\in \M_{\NB,1}$ of gradient flow~\eqref{eq:gradFlowN} is  an SOSP of problem~\eqref{eq:mainN},  for  every balanced initialization $W_{\NB,0}\in \neigh_{\NB,\a}$ outside a subset with  Lebesgue measure zero.

We finally invoke Theorem~\ref{thm:wellBehavedMain}  to conclude that this SOSP $\ol{W}_{\NB}\in \M_{\NB,1}$ is in fact a global minimizer of $L_N$ in $\R^{d_{\NB}}$.
This conclusion is summarized below.

\begin{thm}[Convergence]\label{cor:conjProved}
Gradient flow~\eqref{eq:gradFlowN} converges to a global minimizer of problem~\eqref{eq:mainN} from  every balanced initialization in $\neigh_{\NB,\a} \subset \R^{d_{\NB}}$, outside of a subset with Lebesgue measure zero, see~\eqref{eq:neighWNB}.
Above, $\a\in (\g_Z,1)$, and Assumption~\emph{\ref{assump:key2}} and its notation are in force. 
\end{thm}

A few important remarks are in order. 
{\circled{1}}~Outside the {lazy training} regime  reviewed in Section~\ref{sec:networkDepth}, to our knowledge,  Theorem~\ref{cor:conjProved} is the first convergence result for linear networks, answering the fundamental question of when gradient flow  successfully trains a linear network.  

{\circled{2}}~Indeed, under Assumption~\ref{assump:key2},
Theorem~\ref{cor:conjProved} improves over Theorem~35 in~\cite{bah2019learning} which does \emph{not} guarantee the convergence of gradient flow~\eqref{eq:gradFlowN} to a solution of problem~\eqref{eq:mainN}, and discussed earlier and in  Appendix~\ref{sec:revBahMain}.

{\circled{3}}~This improvement was achieved by provably restricting the gradient flow~\eqref{eq:gradFlowN} to its {stable set} $\neigh_{\NB,\a}$ in~\eqref{eq:neighWNB}, and  such a restriction  is indeed {necessary} as detailed earlier.

{\circled{4}}~{Note that Theorem~\ref{cor:conjProved} sheds light on the theoretical aspects of the training of neural networks, and should not be viewed as an initialization technique for linear networks. In turn, linear networks only serve to improve our theoretical understanding of neural networks in general.}

Let us also examine the content of Assumption~\ref{assump:key2}. {\circled{1}}~The case $r=1$ in~\eqref{eq:rankNet} corresponds to the {spiked covariance} model in statistics, and covers the important case of networks with a scalar output.
Lastly, $r=1$ appears to be the natural building block for extension to $r>1$, which remains an open problem, see the discussion after~\eqref{eq:whitenedData}. 

{\circled{2}}~The assumption of whitened data in~\eqref{eq:whitenedData} is commonly used in the context of linear networks, see for example~\cite{arora2018convergence,bartlett2019gradient}. 
{\circled{3}}~The requirement that $\rank(Z)=\rank(Y\proj_X)\ge r= 1$ in Assumption~\ref{assump:key2}
is evidently necessary to avoid the limit point of zero.

{\circled{4}}~Finally, it is known that the induced flow~\eqref{eq:flowW} for an \emph{unbalanced} initialization
deviates rapidly from its balanced counterpart. 
It is therefore not clear if an unbalanced flow would provably avoid rank-degenerate limit points.
However, we suspect that any disadvantage of an unbalanced initialization will disappear asymptotically as the network depth $N$ grows larger, 
see Equation~8 in~\cite{bah2019learning}.

\section{Convergence Rate of Gradient Flow}\label{sec:networkDepth}

In view of Theorem~\ref{cor:conjProved}, it is natural to ask how fast we can train a linear  network with gradient flow. However, Theorem~\ref{cor:conjProved} is notably silent about the convergence \emph{rate} of gradient flow~\eqref{eq:gradFlowN} to a solution of problem~\eqref{eq:mainN}.
In short, is it  possible for gradient flow to efficiently solve problem~\eqref{eq:mainN}?

As we review now, this fundamental question has not been answered in the literature beyond the lazy training regime.
Indeed, several works have contributed to our understanding here, including \cite{shamir2018exponential,bartlett2019gradient,du2019width,wu2019global,hu2020provable,wu2019global}, and~\cite{gunasekar2017implicit,soudry2018implicit,ji2018gradient,arora2019implicit,rahaman2018spectral,du2018algorithmic} in the related area of {implicit regularization}. 

For our purposes,~\cite{arora2018convergence} exemplifies the current state of the art and its shortcomings. 
Loosely speaking, Theorem~1 in~\cite{arora2018convergence} states that, when the initial loss is \emph{small}, gradient flow~\eqref{eq:gradFlowN} solves problem~\eqref{eq:mainN} to an accuracy of $\epsilon>0$ in the order of
\begin{align}
C^{-(1-\frac{1}{N})} \log(1/\epsilon) \label{eq:aroraRate}
\end{align}
time units, where $C$ is independent of the  depth $N$ of the linear network. For completeness, Theorem~1 in~\cite{arora2018convergence}  is reviewed in Appendix~\ref{sec:lazy}.

Theorem~1 in~\cite{arora2018convergence} might disappoint the researchers.
For one,~\eqref{eq:aroraRate} suggests that increasing the network depth $N$ only marginally speeds up the training.

More concerning is that Theorem~1 in~\cite{arora2018convergence} requires a close initialization, which  is \emph{not} necessary for convergence, see Theorem~\ref{cor:conjProved}. 
Indeed, Theorem~1 in~\cite{arora2018convergence} hinges on a {perturbation argument},  whereby 
the initialization $W_{\NB,0}$ of gradient flow~\eqref{eq:gradFlowN} must satisfy
\begin{align}
    L_N(W_{\NB,0})= \text{sufficiently small} . \qquad \text{(see (\ref{eq:mainN}))}
    \label{eq:goodInitArora}
\end{align}

In this sense,~\cite{arora2018convergence} joins the growing body of literature that quantifies the behavior of neural networks when the \emph{learning trajectory is short}
~\cite{du2018gradient,li2018learning,allen2018convergence,allen2018learning,zou2018stochastic,arora2019fine,tian2017analytical,brutzkus2017globally,
brutzkus2017sgd,du2018power,zhong2017learning,zhang2018learning,wu2019global,shin2019trainability,zhang2019fast,su2019learning,cao2019generalization,chen2019much,oymak2018overparameterized},  to name a few.

To be sure, restricting the initialization  
is   necessary for successful training~\cite{sutskever2013importance}.  For example, gradient flow would    stall when initialized at a  saddle point.

However, it is widely-believed that first-order  algorithms can successfully train neural networks far beyond the {lazy training} regime considered by~\cite{arora2018convergence} and others, and  the line of research exemplified by~\cite{arora2018convergence} is, while  valuable,   highly over-represented in the literature.

Indeed, the learning trajectory of neural networks is in general \emph{not short}, and the learning is often \emph{not local}. We refer to~\cite{chizat2019lazy,yehudai2019power} for a detailed critique of lazy training, see also Appendix~\ref{sec:lazy}.

Let us dub this more general regime \emph{non-local training}. 
Quantifying the non-local convergence rate of linear networks is a vital step towards understanding the non-local training of  neural networks in general.

In an important setting, this section indeed quantifies the non-local training of linear networks, and  addresses both of the shortcomings of Theorem~1 in~\cite{arora2018convergence}.

More specifically, for the case  $r=1$ in~\eqref{eq:rankNet},  Theorem~\ref{thm:MainSimple} below  quantifies the convergence rate of gradient flow~\eqref{eq:gradFlowN} to a solution of problem~\eqref{eq:mainN}, even when~\eqref{eq:goodInitArora} is \emph{violated}.

Moreover, Theorem~\ref{thm:MainSimple} establishes that the faraway convergence rate of gradient flow improves  by increasing the network depth. 
All assumptions for this section are collected in Assumption~\ref{assump:key2}. Let us turn to the details now.

Instead of the convergence rate of gradient flow~\eqref{eq:gradFlowN} to a solution of problem~\eqref{eq:mainN}, we equivalently study the convergence rate of induced flow~\eqref{eq:flowW}, as detailed next. 
The following result is a consequence of Theorem~\ref{cor:conjProved}, proved in Appendix~\ref{sec:proofSwitchFocusRate}.

\begin{lem}[Convergence of induced flow]\label{lem:switchFocusRate}
In the setting of Theorem~\emph{\ref{cor:conjProved}}, if gradient flow~\eqref{eq:gradFlowN} converges to a solution of problem~\eqref{eq:mainN}, then  induced flow~\eqref{eq:flowW} converges to the solution $Z_1$ of problem~\eqref{eq:whyPCA}. Here, $Z_1$ was defined in~\eqref{eq:decomZ0}.
\end{lem}

To quantify the convergence rate of induced flow~\eqref{eq:flowW}, let us define the new loss function
\begin{align}
    L_{1,1}(W) & := \frac{1}{2}\| Z_1 - W\|_F^2.
    \qquad \text{(see~\eqref{eq:decomZ0})}
    \label{eq:newLossDefn}
\end{align}

In this section, we  often opt for subscripts to compactly show the dependence of variables on time $t$, for example, $W_t$ as a shorthand for the induced flow  $W(t)$.  
With $r=1$, recall that induced flow~\eqref{eq:flowW} satisfies $\rank(W_t)\le 1$, see~(\ref{eq:rankNet},\ref{eq:gradFlowToIndFlow}). Assuming that $\rank(W_0)=1$ at  initialization,  the induced flow remains rank-$1$  by Lemma~\ref{lem:rankCteFlow}.
Recall also the analytic SVD of  the induced flow in~\eqref{eq:analyticSVDCVG}. The induced flow thus admits the analytic thin SVD 
\begin{align}
    W_t & \overset{\svd}{=} u_t \cdot s_t \cdot v_t^\top, \qquad \forall t\ge 0,
    \label{eq:analyticTSVD}
\end{align}
where  $u_t\in \R^{d_y}$ and $v_t\in \R^{d_x}$ have unit-norm, and $s_t>0$ is the only nonzero singular value of $W_t$.

A simple calculation using~\eqref{eq:analyticTSVD}, deferred to Appendix~\ref{sec:derLossToRHS}, upper bounds the loss function $L_{1,1}$ in~\eqref{eq:newLossDefn} as
\begin{align}
    L_{1,1}(W_t) & \le  \overset{T_{1,t}}{\overbrace{\frac{1}{2} (s_t - u_t^\top Z_1 v_t)^2}} \nonumber\\
    & \qquad +
    s_Z \underset{T_{2,t}}{\underbrace{(s_Z - u_t^\top Z_1 v_t)}}.
    \label{eq:lossToRHS}
\end{align}
Roughly speaking, $T_{1,t}$ above gauges the error in estimating the (only) nonzero singular value $s_Z$ of the target $Z_1$, whereas $T_{2,t}$ gauges the misalignment between $W_t$ and $Z_1$. Both $T_{1,t},T_{2,t}$ are nonnegative for all $t\ge 0$, see~(\ref{eq:decomZ0},\ref{eq:analyticTSVD}).

To quantify the convergence rate of  induced flow~\eqref{eq:flowW} to the global minimizer $Z_1$ of problem~\eqref{eq:whyPCA}, we next write down the evolution of the loss function $L_{1,1}$ in~\eqref{eq:newLossDefn} as 
\begin{align}
    & \frac{\der L_{1,1}(W_t)}{\der t} \nonumber\\
    & = \l\langle \nabla L_{1,1}(W_t) , \dot{W}_t \r\rangle 
    \qquad \text{(chain rule)} \nonumber\\
    & = - m \l\langle W_t - Z_1, \A_{W_t} ( W_t - Z ) \r\rangle,
    \,\, \text{(see (\ref{eq:flowW},\ref{eq:defnZ}))}
    \label{eq:evolveLossRaw}
\end{align}
where the last line also uses the whitened data in~\eqref{eq:whitenedData}. 
Starting with the definition of $\A_{W_t}$ in~\eqref{eq:defnAW}, we can bound the last line of~\eqref{eq:evolveLossRaw},  see in Appendix~\ref{sec:proofEvolveLemma} for the proof.
\begin{lem}[Evolution of loss]\label{lem:evolveLemma} For induced flow~\eqref{eq:flowW} and the loss function $L_{1,1}$ in~\eqref{eq:newLossDefn}, it 
 holds that 
 \begin{align}
    & \frac{\der L_{1,1}(W_t)}{\der t} \nonumber\\
     & \le  -2m N s_t^{2-\frac{2}{N}} T_{1,t} -2 m s_t^{2-\frac{2}{N}} ( u_t^\top Z_1 v_t) T_{2,t} \nonumber\\
      &  + \sqrt{2} mNs_t^{2-\frac{2}{N}} \g_Z \sqrt{T_{1,t}} T_{2,t}  + 2ms_t^{2-\frac{2}{N}}  s_{Z,2} T_{2,t},
     \label{eq:evolveSimpleLemma}
\end{align}
see~\emph{(\ref{eq:simpler},\ref{eq:analyticTSVD},\ref{eq:lossToRHS})} 
for the notation involved.
\end{lem}

Loosely speaking, the two nonpositive terms on the right-hand side of~\eqref{eq:evolveSimpleLemma} are the contribution of the target matrix~$Z_1$ in~\eqref{eq:decomZ0}, whereas the two nonnegative terms there are the contribution of the residual matrix $Z-Z_1$.  
The (unwanted) nonnegative terms in~\eqref{eq:evolveSimpleLemma} 
 vanish if $Z=Z_1$ is rank-$1$ and, consequently, $\g_Z=s_{Z,2}=0$, see~\eqref{eq:simpler}. 
In view of~(\ref{eq:lossToRHS},\ref{eq:evolveSimpleLemma}), we make two observations: 

    \circled{1}  Both $T_{1,t}$ and $T_{2,t}$ in~\eqref{eq:lossToRHS} appear with negative factors in the dynamics of~\eqref{eq:evolveSimpleLemma}. For loss $L_{1,1}$ to reduce rapidly, we must ensure that $s_t$ and $u_t^\top Z_1 v_t$ both remain bounded away from zero for all $t\ge 0$.
    
    \circled{2}
    $T_{1,t}$ has a large negative factor of $-N$ in the evolution of loss function in~\eqref{eq:evolveSimpleLemma}, and is therefore expected to reduce much faster with time for deeper linear networks.
    
Let us fix $\a\in[\g_Z,1)$ and $\b>1$. Given the first observation above, it is natural to restrict the initialization of gradient flow to a  subset of the feasible set of problem~\eqref{eq:whyPCA}, 
specified as 
\begin{align}
   &  \neigh_{\a,\b}(Z_1):= \Big\{ W \overset{\svd}{=} u_W \cdot s_W \cdot v_W^\top : \nonumber\\
    &\qquad \qquad    (\a-\g_Z) s_{Z}< s_W < \b s_Z, \nonumber\\
    & \qquad \qquad  u_W^\top Z_1 v_W> \a s_{Z}
    \Big\} \subset \R^{d_y\times d_x},
    \label{eq:neighWTwo}
\end{align}
where $s_Z,\g_{Z}$ were defined in~\eqref{eq:simpler}. 

The necessity of such a restriction was discussed after~\eqref{eq:neighWNB}, and the (new) upper bound on $s_W$ in~\eqref{eq:neighWTwo} controls the (unwanted) positive terms in~\eqref{eq:evolveSimpleLemma}. Note that $\neigh_{\a,\b}(Z_1)$ is  a {neighborhood} of $Z_1$, i.e., $Z_1 \in \neigh_{\a,\b}(Z_1)$ 
by~\eqref{eq:decomZ0}.

Once initialized in $\neigh_{\a,\b}(Z_1)$,  induced flow~\eqref{eq:flowW} remains in $\neigh_{\a,\b}(Z_1)$, see Appendix~\ref{sec:proofRemainNeighSimpleTwo},  closely related to Lemma~\ref{lem:remainNeighSimpleNBDomain}. 

\begin{lem}[Stable set]\label{lem:remainNeighSimpleTwo}
Fix $\a\in [\g_Z,1)$ and $\b>1$. For induced flow~\eqref{eq:flowW}, 
$W_0\in \neigh_{\a,\b}(Z_1)$ implies that $W_t\in \neigh_{\a,\b}(Z_1)$ for all $t\ge 0$. 
Above, Assumption~\emph{\ref{assump:key2}} and the notation therein are in force. 
\end{lem}

In view of Lemma~\ref{lem:remainNeighSimpleTwo}, we can now use~(\ref{eq:neighWTwo}) to bound $s_t$ and $u_t^\top Z_1 v_t$  in~\eqref{eq:evolveSimpleLemma}. 
We can then distinguish two regimes  (fast and slow convergence)  in the dynamics of the loss function in~\eqref{eq:evolveSimpleLemma} depending on the dominant term on the right-hand side of~\eqref{eq:lossToRHS}. The remaining technical details 
are  deferred to Appendix~\ref{sec:proofMainSimple} and we finally arrive at the following result.

\begin{thm}[Convergence rate]\label{thm:MainSimple}
    
     With Assumption~\emph{\ref{assump:key2}} and its notation in force, fix $\a\in(\g_Z,1)$ and $\b>1$. Suppose that the inverse spectral gap $\g_Z$ is small enough so that the exponents below are both negative.

     Consider gradient flow~\eqref{eq:gradFlowN} with the balanced initialization $W_{\NB,0}=(W_{1,0},\cdots,W_{N,0})\in \R^{d_{\NB}}$  such that $W_0:=W_{N,0}\cdots W_{1,0}\in \R^{d_y\times d_x}$ satisfies 
     \begin{align}
         & \rank(W_{0}) = 1, \qquad  W_0 \overset{\svd}{=} u_0 s_0 v_0^\top, \nonumber\\
         & (\a-\g_Z) s_Z < s_0 < \b s_Z, \qquad  u_0^\top Z_1 v_0 > \a s_Z. 
         \label{eq:cndInitThmMain}
     \end{align}
    
    Let $W_{\NB}(t)=(W_1(t),\cdots, W_N(t))$ be the output of gradient flow~\eqref{eq:gradFlowN} at time $t$, and set $W(t):=W_N(t)\cdots W_1(t)$, which satisfies $\rank(W(t))=1$ for every $t\ge 0$. 
    
    Let $\tau\ge 0$ be the first time when $s(\tau) \le \sqrt{6} s_Z$, where $s(\tau)$ is the (only) nonzero singular value of $W(\tau)$. Then the distance to the target matrix $Z_1$ in~\eqref{eq:decomZ0} evolves as
\begin{subequations}    
\begin{align}
    &  \forall t \le \tau, \qquad \|Z_1-W(t)\|_F^2  \le \|Z_1-W_0\|_F^2
    \qquad    \label{eq:fastCvgFinal}\\
    &  \qquad \cdot e^{ -m N s_Z ^{2-\frac{2}{N}} \l( (\a-\g_Z)^{2-\frac{2}{N}} - 2\g_Z \b^{2-\frac{2}{N}} \r) t }. \nonumber\\
    &\forall t \ge \tau, \qquad \|Z_1 - W(t)\|_F^2 \le \|Z_1-W(\tau) \|_F^2
          \label{eq:slowCvgFinal}\\
    &  \qquad \cdot e^{ -m s_Z^{2-\frac{2}{N}} 
    \l( \a (\a -\g_Z)^{2-\frac{2}{N}} - 2\g_Z N \b^{2-\frac{2}{N}}\r) (t-\tau) }. \nonumber
\end{align}        
\label{eq:cvgFinalLump}
\end{subequations}

\end{thm}
\vspace{-15pt}
Under Assumption~\ref{assump:key2}, Theorem~\ref{thm:MainSimple} states that gradient flow successfully trains a linear network with linear  rate,  when initialized in the stable set. 

As we will see shortly, Theorem~\ref{thm:MainSimple} is the first result to quantify the convergence rate of gradient flow beyond the widely-studied lazy training regime.
The remarks after Theorem~\ref{cor:conjProved} again apply here about Assumption~\ref{assump:key2} and  the case~$r>1$. A few  additional remarks are in order. 
 
{\circled{1}} Rephrasing (\ref{eq:cvgFinalLump}), gradient flow~\eqref{eq:gradFlowN} solves problem~\eqref{eq:mainN} to an accuracy of $\epsilon>0$ in the order of 
\begin{align*}
    \begin{cases}
    \frac{1}{mNs_Z^2}\l(   (\a-\g_Z)^{2} - 2\g_Z \b^{2} \r)^{-1} \log(C/\epsilon) & \epsilon > \epsilon_0 \\
    \frac{1}{ms_Z^2}\l(   \a (\a -\g_Z)^{2} - 2\g_Z N \b^{2} \r)^{-1} \log(C/\epsilon)\\
    - \tau\l(N \frac{   (\a-\g_Z)^{2} - 2\g_Z \b^{2} }{\a (\a -\g_Z)^{2} - 2\g_Z N \b^{2}}  -1\r)  &  \epsilon\le \epsilon_0
    \end{cases}
\end{align*}
time units.  Above, $\epsilon_0 $ is the right-hand side of~\eqref{eq:slowCvgFinal}, evaluated at $t=\tau$.

{\circled{2}}~In Theorem~\ref{thm:MainSimple}, the end-to-end initialization matrix $W_0$ in~\eqref{eq:cndInitThmMain} is positively correlated with the target matrix $Z_1$, and  away from the origin, see our earlier discussions for the necessity of such restricted initialization.
Note also that \eqref{eq:cndInitThmMain} should not be seen as an initialization scheme but as a theoretical result.

{\circled{3}}~The faraway convergence rate in~\eqref{eq:cvgFinalLump} improves with increasing the network depth $N$, whereas the nearby convergence rate does not appear to benefit from increasing $N$. This improved faraway convergence rate should be contrasted with Arora's result in~\eqref{eq:aroraRate}.



{\circled{4}}~Crucially, the lazy training results fail to apply here. To see this, with the initialization $W_{\NB,0}=(W_{1,0},\cdots,W_{N,0})$, Claim~1 in~\cite{arora2018convergence} uses a perturbation argument,
which requires that 
\begin{align}
    \| Z - W_{N,0}\cdots W_{1,0}\|_F < s_{\min}(Z), 
\label{eq:flaw}
\end{align}
which is \emph{impossible} unless trivially $\rank(Z)\le r$. Indeed, the network architecture  forces that $\rank(W_{N,0}\cdots W_{1,0})\le r$, see~\eqref{eq:rankNet}. 

In contrast, Theorem~\ref{thm:MainSimple} applies even when~\eqref{eq:flaw} is violated, as it does away entirely with the limitations of a perturbation argument.
Theorem~\ref{thm:MainSimple}  thus ventures beyond the reach of the {lazy training} regime in~\cite{arora2018convergence}, which has dominated the recent literature of neural networks, thus signifying the importance of this breakthrough.

\begin{figure}[h!]
   \centering
     \includegraphics[scale=.45]{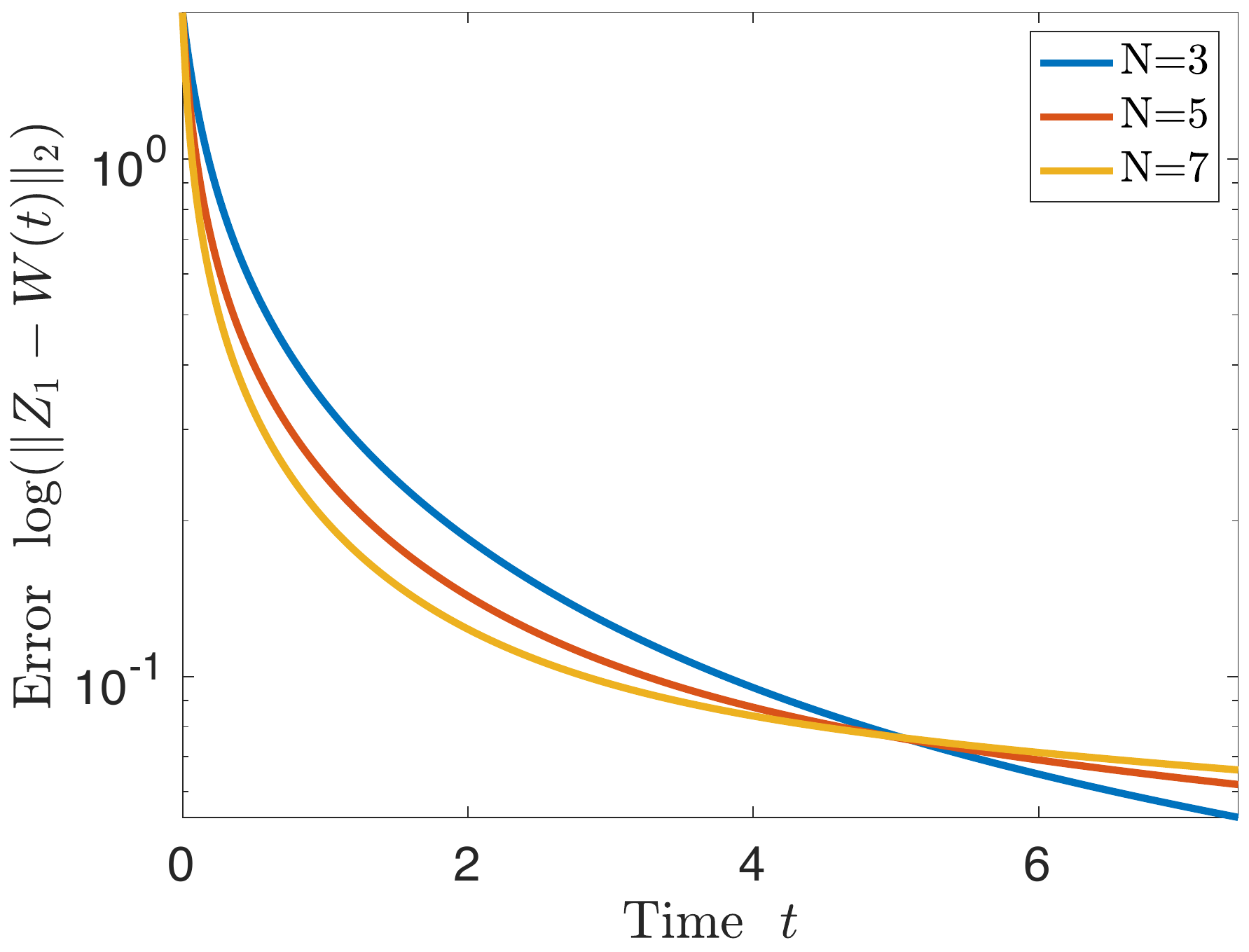}
    \caption{
    Suppose that the sample size is $m=50$, and  consider a randomly-generated  whitened training dataset 
    $
    (X,Y)\in \R^{d_x\times m}\times \R^{d_y\times m},
    $
    with $d_x=5$ 
    and $d_y=1$. For this dataset, the  above figure depicts the distance from induced flow~\eqref{eq:flowW} to the target vector $Z_1=Z=YX^\top/m$ in~(\ref{eq:defnZ},\ref{eq:decomZ0}), plotted versus time $t$, for training a linear network with $d_x$ inputs and $d_y$ outputs, as the network depth $N$ varies. 
    \vspace{5pt}
    \\
    The direction of the  initial end-to-end  vector $W_0\in \R^{d_y\times d_x}$ is obtained by randomly rotating the direction of the  target vector $Z_1$ by about $30$ degrees. We also set $\|W_0\|_2=10\|Z\|_2$. 
    Instead of induced flow~\eqref{eq:flowW}, we implemented  the discretization of~\eqref{eq:flowW} obtained from the explicit (or forward) Euler method with a step size of $10^{-6}$ with $10^5$ steps.  
    \vspace{5pt}
    \\
    This simple numerical example  visualizes the (gradual) slow-down in the convergence rate of gradient flow with time, see~(\ref{eq:cvgFinalLump}), and also shows the faster faraway convergence rate for deeper networks, see Theorem~\ref{thm:MainSimple}. 
    The above figure also suggests that the nearby convergence rate of gradient flow~\eqref{eq:gradFlowN} might actually be slower for deeper networks. It is however difficult to theoretically infer this    from Theorem~\ref{thm:MainSimple}, because~\eqref{eq:slowCvgFinal} is an \emph{upper bound} for the nearby  error. 
      \vspace{5pt}
    \\
    The precise nearby convergence rates of linear networks (and any trade-offs associated with the network depth) thus remain as open questions. Note also that the local analysis of~\cite{arora2018convergence} cannot be applied here, as discussed after Theorem~\ref{thm:MainSimple}.   
    }
     \label{fig:demo}
 \end{figure}

Thorough numerics for linear networks are abound, see for example~\cite{bah2019learning,arora2018convergence,arora2018optimization}, and we refrain from lengthy simulations  and only provide a numerical example in Figure~\ref{fig:demo}.
to visualize the (gradual) change of regimes from fast to slow convergence, see~(\ref{eq:fastCvgFinal},\ref{eq:slowCvgFinal}). This example also suggests new research questions about linear networks.  

\newpage
\twocolumn


\newpage
\appendix

\section{Derivation of~(\ref{eq:mainOne},\ref{eq:mainTwo})}\label{sec:surjective}

To show~(\ref{eq:mainOne},\ref{eq:mainTwo}), it suffices to show that the map 
\begin{align}
    & \map(d_{\NB}) : \R^{d_{\NB}} \rightarrow \M^{d_N\times d_0}_{1,\cdots,r}
    \nonumber\\
    & W_{\NB}=(W_1,\cdots,W_N) \rightarrow W=W_N\cdots W_1,
    \label{eq:mapDefn}
\end{align}
is surjective, which we now set out to do. 

Above, $d_{\NB}=(d_0,\cdots,d_N)$ and  $\R^{d_{\NB}}=\R^{d_0}\times \cdots \times \R^{d_N}$ is the domain of the function. Also, $\M^{d_N\times d_0}_{1,\cdots,r} \subset \R^{d_N\times d_0}$ is the set of all $d_N\times d_0$ matrices of rank at most $r$. As a side note, $\M^{d_N\times d_0}_{1,\cdots,r}$ is the closure of the manifold of rank-$r$ matrices. Lastly, the network architecture dictates that
satisfies 
\begin{align}
\min_j d_j = r. \qquad \text{(see~\eqref{eq:rankNet})}
\label{eq:rankNetRe}
\end{align}

The proof of this surjective property is by induction. 

The base of induction for $N=1$ is trivial because $\map(d_0,d_1)$ is simply the identity map  by~\eqref{eq:mapDefn} and thus  surjective, in particular for any pair of integers $(d_0,d_1)$ that satisfies~\eqref{eq:rankNetRe}.

For the step of induction, suppose that $\map(d_{\NB})$ is surjective for every tuple $d_{\NB}=(d_0,\cdots,d_N)$ that satisfies~\eqref{eq:rankNetRe}.

For an arbitrary integer $d_{N+1}$, consider also an arbitrary matrix 
\begin{align}
W\in \M_{1,\cdots,r}^{d_{N+1}\times d_0},
\label{eq:newW}
\end{align}
with the  SVD
\begin{align}
    W \overset{\text{SVD}}{=}\wt{U} \cdot \wt{S}\wt{V}^\top =: \wt{U} \cdot \wt{Q},
    \label{eq:decompWproof}
\end{align}
where $\wt{U}\in \R^{d_{N+1}\times d_{N+1}}$ and $\wt{V}\in \R^{d_0\times d_0}$ are orthonormal bases, and $\wt{S}\in\R^{d_{N+1}\times d_0}$ contains the singular values of~$W$. 

In particular, note that $\wt{Q}\in \R^{d_{N+1}\times d_0}$. 
Note also that \eqref{eq:newW} implies that  
\begin{align}
\rank(\wt{Q})=\rank(W) \le r,
\label{eq:QWr}
\end{align}
because $\wt{U}$ is an orthonormal basis. 

Combining~\eqref{eq:rankNetRe} and~\eqref{eq:QWr}, we reach
\begin{align}
    \rank(\wt{Q}) \le r \le d_N.\label{eq:rSandwitched}
\end{align}

In view of~\eqref{eq:rSandwitched}, it is therefore possible (by padding with zero columns or removing some columns from $\wt{U}$ and the corresponding rows from $\wt{Q}$) to create $\wt{U}'$ and $\wt{Q}'$ such that 
\begin{align}
     & W  = \wt{U}'\cdot \wt{Q}'^\top, \nonumber\\
    & \text{where}\,\, \wt{U}'\in \R^{d_{N+1}\times d_N}, \,\, \wt{Q}'\in \R^{d_N\times d_0}. 
    \label{eq:decompWproof2}
\end{align}

In this construction, $\rank(\wt{Q}')=\rank(\wt{Q})\le r$ and $\wt{Q}'\in \R^{d_N\times d_0}$. Consequently, the step of induction guarantees the existence of $W_{\NB}=(W_N,\cdots,W_1)\in \R^{d_{\NB}}$ such that 
\begin{align}
    \wt{Q}' = W_N \cdots W_1. \label{eq:stepUsed}
\end{align}
It  follows that 
\begin{align}
    W & = \wt{U}'\cdot \wt{Q}'  \qquad \text{(see \eqref{eq:decompWproof2})} \nonumber\\
    & = \wt{U}' W_N \cdots W_1. \qquad \text{(see \eqref{eq:stepUsed})}
\end{align}
That is, 
\begin{align}
    W = \map(d_0,\cdots,d_{N+1})[W_1,\cdots,W_N,\wt{U}'],
\end{align}
which completes the induction. We thus proved that $\Pi(d_{\NB})$ is a surjective map for every tuple $d_{\NB}$ that satisfies~\eqref{eq:rankNetRe}.

\section{Proof of Lemma~\ref{lem:LNtoL2}}\label{sec:proofLNtoL2}

Let $\ol{W}_{\NB}=(\ol{W}_1,\cdots,\ol{W}_N)\in \R^{d_{\NB}}$ be an FOSP of problem~\eqref{eq:mainN}. For an infinitesimally small perturbation $\D_{\NB}=(\D_1,\cdots, \D_N)\in \R^{d_{\NB}}$, we can expand $L_N$ in~\eqref{eq:mainN} as
\begin{align}
    & L_N(\ol{W}_{\NB}+\D_{\NB}) \nonumber\\
    & = 
    L_N(\ol{W}_{\NB})+ \nabla L_N(\ol{W}_{\NB})[\D_{\NB}] \nonumber\\
    & + \frac{1}{2} \nabla^2 L_N(\ol{W}_{\NB})[\D_{\NB}]+o \nonumber\\
    & = L_N(\ol{W}_{\NB})+\frac{1}{2}\nabla^2 L_N(\ol{W}_{\NB})[\D_{\NB}]+o.
    \label{eq:firstExpand}
\end{align}
where $o$ represents (negligible) higher order terms, and the second identity above holds because   $\ol{W}_{\NB}$ is assumed to be an FOSP  in Lemma~\ref{lem:LNtoL2}, see Definition~\ref{defn:fosp}. Above, $\nabla^2 L_N(\ol{W}_{\NB})[\D_{\NB}]$ contains all second order terms in the variables $\D_{\NB}$. 

Let $j_0$  correspond to a layer  with the smallest width within the linear network~(\ref{eq:linNet},\ref{eq:linearNetMap}), i.e.,
\begin{align}
    r & = \min_{j\le N} d_j  \qquad \text{(see \eqref{eq:rankNet})} \nonumber\\
    & = d_{j_0}. 
\end{align}
We also set 
\begin{align}
   & \ol{P} := \ol{W}_N \cdots \ol{W}_{j_0+1} \in \R^{d_y\times r},\nonumber\\
   & \ol{Q} := \ol{W}_{j_0} \cdots \ol{W}_1\in \R^{r\times d_x},
\end{align}
for short, and  note that 
\begin{align}
    & \ol{W}  := \ol{P}  \cdot \ol{Q} = \ol{W}_N \cdots \ol{W}_1,\nonumber\\
    & \rank(\ol{W}) = \rank(\ol{P}) = \rank(\ol{Q}) = r,
    \label{eq:PQrankR}
\end{align}
where the second line above holds by  the assumption of Lemma~\ref{lem:LNtoL2}. 

Indeed, $\ol{W}=\ol{P}\cdot \ol{Q}$ implies that
\begin{align}
\min(\rank(\ol{P}),\rank(\ol{Q}))\ge r.
\label{eq:uppBndMinRank}
\end{align}
Note also that $\ol{P}$ has $r$ columns and $\ol{Q}$ has $r$ rows, thus 
\begin{align}
\max(\rank(\ol{P}),\rank(\ol{Q})) \le r.
\label{eq:lowBndMaxRank}
\end{align}
Together,~(\ref{eq:uppBndMinRank}) and~(\ref{eq:lowBndMaxRank}) give the second line of~\eqref{eq:PQrankR}.

On the one hand, for an arbitrary $(\D_P,\D_Q)$, we can relate the perturbation of $\ol{W}_{\NB}$ to the perturbation of $(\ol{P},\ol{Q})$ as 
\begin{align}
    & (\ol{W}_N+\D_N)\cdots (\ol{W}_1+\D_1) \nonumber\\
    & = (\ol{P}+\D_P)(\ol{Q}+\D_Q),
    \label{eq:prodsMatch}
\end{align}
where 
\begin{align}
    & \D_1 = \ol{W}_1\ol{Q}^\dagger \D_Q,\nonumber\\
    & \D_i = 0, \qquad  2\le i\le N-1,\nonumber\\
    & \D_N = \D_P \ol{P}^\dagger \ol{W}_N, 
    \label{eq:mapProof}
\end{align}
and $\dagger$ denotes the pseudo-inverse, and we used the  second identity in~\eqref{eq:PQrankR}. 

Indeed, for the choice of $\D_{\NB}$ in~\eqref{eq:mapProof}, it holds that 
\begin{align}
    & (\ol{W}_N+\D_N) \cdots (\ol{W}_1+\D_1) \nonumber\\
    & = (\ol{W}_N+\D_P \ol{P}^\dagger \ol{W}_N) \ol{W}_{N-1} \cdots \nonumber\\
    & \qquad \cdots  \ol{W}_2 (\ol{W}_1+\ol{W}_1\ol{Q}^\dagger \D_Q) \nonumber
    \qquad \text{(see \eqref{eq:mapProof})}
    \\
    & = (I_{d_y} + \D_P \ol{P}^\dagger ) \ol{W}_N \ol{W}_{N-1}\cdots \nonumber\\
    & \qquad \cdots \ol{W}_2 \ol{W}_1 (I_{d_x}+ \ol{Q}^\dagger \D_Q)\nonumber\\
    & = (I_{d_y} + \D_P \ol{P}^\dagger ) \ol{P} \cdot \ol{Q} (I_{d_x}+ \ol{Q}^\dagger \D_Q) \quad \text{(see \eqref{eq:PQrankR})} \nonumber\\
    & = (\ol{P}+\D_P) (\ol{Q}+\D_Q),
\end{align}
which agrees  with~\eqref{eq:prodsMatch}. The last line above uses the second identity in~\eqref{eq:PQrankR}, i.e., $\rank(\ol{P})=\rank(\ol{Q})=r$. Above, $I_{d_y}\in \R^{d_y\times d_y}$ is the identity matrix. 

On the other hand, we can  expand $L_2$ in~\eqref{eq:mainOne} as 
\begin{align}
    & L_2(\ol{P}+\D_P,\ol{Q}+\D_Q) \nonumber\\
    & = L_2(\ol{P},\ol{Q})+ \nabla L_2 (\ol{P},\ol{Q})[\D_P,\D_Q] \nonumber\\
    & + \nabla^2  L_2 (\ol{P},\ol{Q})[\D_P,\D_Q]+o,
    \label{eq:secondExpand}
\end{align}
where $o$ again higher order terms. Above, $\nabla L_2 (\ol{P},\ol{Q})[\D_P,\D_Q] $ collects all first order terms in the variables $(\D_P,\D_Q)$. Likewise, $\nabla^2 L_2 (\ol{P},\ol{Q})[\D_P,\D_Q] $ contains all second order terms in $(\D_P,\D_Q)$.

For convenience, let us  define the map
\begin{align}
    L_1: & \R^{d_y\times d_x} \rightarrow \R \nonumber\\
    & W \rightarrow \frac{1}{2}\|Y-WX\|_F^2,
\end{align}
and note that 
\begin{align}
    & L_2(P,Q) = L_1(PQ), \qquad \text{(see \eqref{eq:mainTwo})} \nonumber\\
    & L_N(W_{\NB}) = L_1(W_N\cdots W_1), \qquad \text{(see \eqref{eq:mainN})}
    \label{eq:defnL2pLNp}
\end{align}
for every $P,Q,W_{\NB}$.

In view of~\eqref{eq:defnL2pLNp}, we now write that 
\begin{align}
    & L_2(\ol{P}+\D_P,\ol{Q}+\D_Q) 
    \nonumber\\
    & = L_1((\ol{P}+\D_P)(\ol{Q}+\D_Q))\qquad \text{(see \eqref{eq:defnL2pLNp})} \nonumber\\
    & = L_1( (\ol{W}_N+\D_N) \cdots (\ol{W}_1+\D_1) ) 
    \qquad \text{(see \eqref{eq:prodsMatch})} 
    \nonumber\\
    & = L_N(\ol{W}_{\NB}+\D_{\NB}),\qquad \text{(see \eqref{eq:defnL2pLNp})}
    \label{eq:preTaylorMatch}
\end{align}
for $\D_{\NB}=(\D_1,\cdots,\D_N)$ specified in~\eqref{eq:mapProof}.

As a result of~\eqref{eq:preTaylorMatch}, the expansions in~\eqref{eq:firstExpand} and~\eqref{eq:secondExpand} must match. That is, for an arbitrary $(\D_{P},\D_Q)$ and the corresponding choice of $\D_{\NB}$ in~\eqref{eq:prodsMatch}, it holds that 
\begin{align}
    & \nabla L_2(\ol{P},\ol{Q})[\D_P,\D_Q]  \nonumber\\
    & = \nabla L_N(\ol{W}_{\NB})[\D_{\NB}] = 0, \qquad \text{(see (\ref{eq:firstExpand},\ref{eq:secondExpand}))}
    \label{eq:match1}
\end{align}
and
\begin{align}
   &  \nabla^2 L_2(\ol{P},\ol{Q})[\D_P,\D_Q] \nonumber\\
   & = \nabla^2 L_N(\ol{W}_{\NB})[\D_{\NB}].\qquad \text{(see (\ref{eq:firstExpand},\ref{eq:secondExpand}))}
   \label{eq:match2}
\end{align}
It follows from~(\ref{eq:match1}) that $(\ol{P},\ol{Q})$ is an FOSP of problem~\eqref{eq:mainOne} if $\ol{W}_{\NB}$ is an FOSP of problem~\eqref{eq:mainN}.

Moreover, if $\ol{W}_{\NB}$ is an SOSP of problem~\eqref{eq:mainN}, then the last line of~\eqref{eq:match2} is nonnegative, see Definition~\ref{defn:secondStN}. That is, 
\begin{align}
   &  \nabla^2 L_2(\ol{P},\ol{Q})[\D_P,\D_Q] \nonumber\\
   & = \nabla^2 L_N(\ol{W}_{\NB})(\D_{\NB})\ge 0,
   \label{eq:match3}
\end{align}
for an arbitrary $(\D_P,\D_Q)$ and the corresponding choice of $\D_{\NB}$ in~\eqref{eq:prodsMatch}. Therefore, $(\ol{P},\ol{Q})$ is an SOSP of problem~\eqref{eq:mainOne} if $\ol{W}_{\NB}$ is an SOSP of problem~\eqref{eq:mainN}. This completes the proof of Lemma~\ref{lem:LNtoL2}.

\section{Proof of Lemma~\ref{lem:L2toPCA}}\label{sec:proofL2toPCA}

Recall that $\proj_X$ and $\proj_{X^\perp}$ denote the orthogonal projections onto the row span of $X$ and its complement, respectively.

Using the decomposition $Q'=Q'\proj_X+Q' \proj_{X^\perp}$, the last program in~\eqref{eq:pca} can be written as 
\begin{align}
    & \underset{P,Q'}{\min} \,\,\frac{1}{2}\|Y\proj_X - PQ'\|_F^2 \nonumber\\
    & = \underset{P,Q'}{\min}\,\, \frac{1}{2}\|Y\proj_X - PQ'\proj_X\|_F^2 + \frac{1}{2}\| PQ'\proj_{X^\perp}\|_F^2.
    \label{eq:preEvident}
\end{align}
From the above decomposition, it is evident that the minimum above is achieved when the last term in~\eqref{eq:preEvident} vanishes. This observation allows us to write that 
\begin{align}
    & \underset{P,Q'}{\min} \,\, \frac{1}{2}\|Y\proj_X - PQ'\|_F^2 \nonumber\\
    & = \underset{P,Q'}{\min}\,\, \frac{1}{2}\|Y\proj_X - PQ'\proj_X\|_F^2 
    \qquad \text{(see \eqref{eq:preEvident})}
    \nonumber\\
    & = \begin{cases}
    \underset{P,Q'}{\min}\,\, \frac{1}{2}\|Y\proj_X - PQ''\|_F^2\\
    \text{subject to }Q'' = Q' \proj_X
    \end{cases} \nonumber\\
    & = \begin{cases}
    \underset{P,Q''}{\min}\,\, \frac{1}{2}\|Y\proj_X - PQ''\|_F^2\\
    \text{subject to }\rowspan(Q'') \subseteq \rowspan(X)
    \end{cases} \nonumber\\
    & = \begin{cases}
    \underset{P,Q''}{\min}\,\, \frac{1}{2}\|Y\proj_X - PQ''\|_F^2\\
    \text{subject to } Q'' = Q X
    \end{cases} \nonumber\\
    & = \underset{P,Q}{\min}\,\, \frac{1}{2}\|Y\proj_X - PQX\|_F^2,
    \label{eq:Evident}
\end{align}
which proves the tight relaxation claimed in~\eqref{eq:pca}. The third identity above uses the fact that the map 
\begin{align}
    & \R^{r\times m} \rightarrow \rowspan(X) \nonumber\\
    & Q'\rightarrow Q'' = Q'\proj_X 
\end{align}
is surjective.

To prove the second claim in Lemma~\ref{lem:L2toPCA}, let $(\ol{P},\ol{Q})$ be an FOSP of problem~\eqref{eq:mainTwo}, which satisfies
\begin{align}
    & 0 = (Y-\ol{P}\cdot \ol{Q}X)X^\top \ol{Q}^\top , \nonumber\\
    & 0 = \ol{P}^\top (Y-\ol{P}\cdot \ol{Q}X) X^\top .
    \label{eq:aboveIdUsed}
\end{align}
After setting 
\begin{align}
    \ol{Q}'=\ol{Q}X,
    \label{eq:QtoQpProof}
\end{align} 
the above identities read as
\begin{align}
    0 & =(Y-\ol{P}\cdot \ol{Q}X) X^\top \ol{Q}^\top  
    \qquad \text{(see \eqref{eq:aboveIdUsed})}
    \nonumber\\
    & = (Y\proj_X - \ol{P}\cdot \ol{Q} X) X^\top \ol{Q}^\top
    \nonumber\\
    & = (Y\proj_X - \ol{P}\cdot \ol{Q}') \ol{Q}'^\top,
    \label{eq:firstIdPr}
\end{align}
and
\begin{align}
  0 &= \ol{P}^\top (Y-\ol{P}\cdot \ol{Q} X)X^\top 
  \qquad \text{(see \eqref{eq:aboveIdUsed})}
  \nonumber\\
  & = \ol{P}^\top (Y \proj_X -\ol{P}\cdot \ol{Q} X)X^\top \nonumber\\
  & = \ol{P}^\top (Y \proj_X -\ol{P}\cdot \ol{Q}') X^\top. 
      \label{eq:secondIdPr1}
\end{align}
Recall that 
\begin{align}
\rowspan(\ol{Q}')\subseteq \rowspan(X).
\qquad \text{(see \eqref{eq:QtoQpProof})}
\label{eq:QpSpanX}
\end{align}
With this in mind,~\eqref{eq:secondIdPr1} implies that 
\begin{align}
    0 & = \ol{P}^\top (Y\proj_X -\ol{P}\cdot \ol{Q}') X^\top 
    \qquad \text{(see \eqref{eq:secondIdPr1})}\nonumber\\
    & = \ol{P}^\top (Y\proj_X -\ol{P}\cdot \ol{Q}'),
    \qquad \text{(see \eqref{eq:QpSpanX})}
          \label{eq:secondIdPr}
\end{align}
where we also used the assumption that $XX^\top$ is invertible.

By combining~(\ref{eq:firstIdPr},\ref{eq:secondIdPr}), we conclude that $(\ol{P},\ol{Q}')$ is an FOSP of problem~\eqref{eq:pca} if $(\ol{P},\ol{Q})$ is an FOSP of problem~\eqref{eq:mainTwo}. 

To prove the last claim of Lemma~\ref{lem:L2toPCA}, let $(\ol{P},\ol{Q})$ be an SOSP of problem~\eqref{eq:mainTwo}, which satisfies 
\begin{align}
    & \frac{1}{2}\|\D_P \ol{Q}X + \ol{P} \D_Q X \|_F^2   + \langle \ol{P}\cdot \ol{Q}X - Y, \D_P\D_Q X \rangle =  \nonumber\\
    &  \frac{1}{2}\|\D_P \ol{Q}X + \ol{P} \D_Q X \|_F^2   + \langle \ol{P}\cdot \ol{Q}X - Y\proj_X , \D_P\D_Q X \rangle \nonumber\\ 
    & \ge 0
    , 
      \qquad \forall (\D_P,\D_Q).
    \label{eq:tempPrSecondOrder}
\end{align}
Let us set $\ol{Q}'=\ol{Q}X$ as before, and also note that the map
\begin{align}
    & \R^{r\times d_x} \rightarrow \rowspan(X) \nonumber\\
    & \D_Q \rightarrow \D_{Q'} = \D_Q X
\end{align}
is evidently surjective. 
Then we may rewrite~\eqref{eq:tempPrSecondOrder} as
\begin{align}
        & \frac{1}{2}\|\D_P \ol{Q}' + \ol{P} \D_{Q'} \|_F^2 + \langle \ol{P}\cdot \ol{Q}' - Y\proj_X , \D_P \D_{Q'}\rangle \ge 0, \nonumber\\
     & \qquad \forall (\D_P,\D_{Q'}) \in \R^{d_y\times r} \times \rowspan(X),
     \label{eq:finalPCAPr}
\end{align}
On the other hand, recall again~\eqref{eq:QpSpanX}. When 
\begin{align}
\D_{Q'} \perp  \rowspan(X),
\label{eq:orthPart}
\end{align}
we have that 
\begin{align}
    & \frac{1}{2}\| \D_P \ol{Q}' + \ol{P} \D_{Q'}\|_F^2  + 
    \langle \ol{P}\cdot \ol{Q}' - Y \proj_X, \D_P \D_{Q'}\rangle
    \nonumber\\
    & = \frac{1}{2} \| \D_P \ol{Q}' \|_F^2 + \|\ol{P} \D_{Q'}\|_F^2 \ge 0, \nonumber\\
    & \forall   (\D_P,\D_{Q'}) \in \R^{d_y\times r} \times \rowspan(X)^\perp ,
         \label{eq:finalPCAPr2}
\end{align}
where the identity above uses~(\ref{eq:QpSpanX},\ref{eq:orthPart}). By combining~(\ref{eq:finalPCAPr},\ref{eq:finalPCAPr2}), we reach
\begin{align}
        & \frac{1}{2}\| \D_P \ol{Q}' + \ol{P} \D_{Q'}\|_F^2  + 
    \langle \ol{P}\cdot \ol{Q}' - Y \proj_X, \D_P \D_{Q'}\rangle \nonumber\\
    & \ge 0, \qquad  \forall (\D_P,\D_{Q'}).
    \label{eq:finalPCAPr3}
\end{align}
It is evident from~(\ref{eq:finalPCAPr3}) that  $(\ol{P},\ol{Q}')$ is an SOSP of problem~\eqref{eq:pca} if $(\ol{P},\ol{Q})$ is an SOSP of problem~\eqref{eq:mainTwo}. This completes the proof of Lemma~\ref{lem:L2toPCA}. 

\section{Proof of Theorem \ref{thm:wellBehavedMain}}\label{sec:proofThmMain}

We begin with a  technical lemma below, proved with the aid of EYM Theorem~\ref{thm:eym}. This result is standard but a proof is included for completeness. 
\begin{lem}\label{lem:noRankDegStPt}
If $\rank(Y\proj_X)\ge r$, then any SOSP $(\ol{P},\ol{Q}')$ of problem~\eqref{eq:pca}  is a global minimizer of problem~\eqref{eq:pca} and satisfies  
\begin{align}
 \rank(\ol{P})=\rank(\ol{Q}')= \rank(\ol{W}) = r,
 \label{eq:rankNotDegPrThm}
\end{align}
where $\ol{W} = \ol{P}\cdot \ol{Q}'$. 
\end{lem}
Before proving the above lemma in the next appendix, let us show how it can be used to prove Theorem~\ref{thm:wellBehavedMain}.

Let us assume that $\rank(Y\proj_X)\ge r$, so that Lemma~\ref{lem:noRankDegStPt} is in force. Then any SOSP $(\ol{P},\ol{Q}')$ of problem~\eqref{eq:pca} is  a global minimizer of problem~\eqref{eq:pca} and satisfies~\eqref{eq:rankNotDegPrThm}.

Let us also assume that $XX^\top$ is invertible, so that Lemma \ref{lem:cmobinedLand} is in force.
Lemma~\ref{lem:cmobinedLand} then implies that any SOSP $\ol{W}_{\NB}$  of problem~\eqref{eq:mainN} corresponds to an SOSP $(\ol{P},\ol{Q}')$ of problem~\eqref{eq:pca}, 
provided that $\ol{W}=\ol{W}_N\cdots\ol{W}_1$ is rank-$r$. 
The relationship between these quantities is
\begin{align}
    \ol{W}_N \cdots \ol{W}_1 X = \ol{W} X = \ol{P} \cdot \ol{Q}'.
    \qquad \text{(see (\ref{eq:PQrankR},\ref{eq:QtoQpProof})}
    \label{eq:relationShipMainProof}
\end{align}

In light of the preceding paragraph, we observe that any SOSP $\ol{W}_{\NB}$ of problem~\eqref{eq:pca} corresponds to a global minimizer $(\ol{P},\ol{Q}')$ of problem~\eqref{eq:pca}, provided that $\ol{W}$ is rank-$r$.

Using the decomposition $Y=Y\proj_X+Y\proj_{X^\perp}$, we can therefore write that 
\begin{align}
    & \frac{1}{2}\| Y - \ol{W}_N\cdots \ol{W}_1 X \|_F^2  \nonumber\\
    &= \frac{1}{2}\| Y\proj_X - \ol{W}_N\cdots \ol{W}_1 X \|_F^2 + \frac{1}{2}\|Y \proj_{X^\perp}\|_F^2 \nonumber\\
    & = \frac{1}{2}\| Y \proj_X - \ol{P} \cdot \ol{Q}'\|_F^2 +  \frac{1}{2}\|Y\proj_{X^\perp}\|_F^2 
    \qquad \text{(see \eqref{eq:relationShipMainProof})}
    \nonumber\\
    & = \underset{P,Q'}{\min} \,\, \frac{1}{2}\| Y\proj_X - PQ'\|_F^2 +  \frac{1}{2}\|Y\proj_{X^\perp}\|_F^2 \nonumber\\
     & = \min_{P,Q} \,\, \frac{1}{2}\| Y - PQX\|_F^2  
     \qquad \text{(see \eqref{eq:pca})}
     \nonumber\\
     & = \underset{W_1,\cdots,W_N}{\min} \,\, \frac{1}{2}\|Y- W_N\cdots W_1 X \|_F^2 . \,\,\, \text{(see \eqref{eq:mainTwo})}
\end{align} 
That is, any SOSP $\ol{W}_{\NB}$  of problem~\eqref{eq:mainN} is a global minimizer of problem~\eqref{eq:mainN}, provided that $\ol{W}$ is rank-$r$. 
This completes the proof of Theorem~\ref{thm:wellBehavedMain}.

\subsection{Proof of Lemma \ref{lem:noRankDegStPt}} \label{sec:proofNoRankDegStPt}

We conveniently assume that  
\begin{align}
    \rank(Y\proj_X) = r, \label{eq:notRankDegRe}
\end{align}
but the same argument is valid also when $\rank(Y\proj_X) > r$. 
Let 
\begin{align}
Y\proj_X \overset{\svd}{=} \wt{U}\wt{S}\wt{V}
\label{eq:svdYPX}
\end{align}
denote the thin SVD of $Y\proj_X$, where
$\wt{U} \in \R^{d_y\times r}$ has orthonormal columns,    $\wt{V}\in \R^{r \times m}$ has orthonormal rows, and the diagonal matrix 
$ 
\wt{S}\in \R^{r\times r}
$ 
contains the singular values of $Y\proj_X$.

By the way of  contradiction, suppose that $(\ol{P},\ol{Q}')$ is an SOSP of problem~\eqref{eq:pca} such that  
\begin{align}
\rank(\ol{P}\cdot \ol{Q}') < r.\label{eq:counterAssumption}
\end{align}
Without loss of generality, we can in fact replace~\eqref{eq:counterAssumption} with 
\begin{align}
    \rank(\ol{P}) =\rank(\ol{Q}') = \rank(\ol{P}\cdot \ol{Q}') < r.
    \label{eq:counterAssumption2}
\end{align}
(Indeed, for example if $\rank(\ol{P})< \rank(\ol{Q}') < r$, then $(\ol{P} \proj_{\ol{S}} ,\proj_{\ol{S}} \ol{Q}')$ takes the same objective value in problem~\eqref{eq:pca} as $(\ol{P},\ol{Q}')$. Here,  $\proj_{\ol{S}}$ is the orthogonal projection onto the subspace $\ol{S}=\rowspan(\ol{P})\cap\text{column span}(\ol{Q}')$. 

On the other hand, by EYM Theorem~\ref{thm:eym}, the SOSP $(\ol{P},\ol{Q}')$ is in fact a global minimizer of problem~\eqref{eq:pca}. Therefore, $( \ol{P} \proj_{\ol{S}},\proj_{\ol{S}}\ol{Q}'$) too is a global minimizer of problem~\eqref{eq:pca} and a fortiori an SOSP of problem~\eqref{eq:pca}. We can thus replace $(\ol{P},\ol{Q}')$ with the SOSP $( \ol{P}\proj_{\ol{S}},\proj_{\ol{S}}\ol{Q}')$ which satisfies $\rank(\ol{P}\proj_{\ol{S}}) = \rank(\proj_{\ol{S}} \ol{Q}')<r$. That is, the assumption made in~\eqref{eq:counterAssumption2} indeed does not reduce the generality of the following argument.)

Assuming~\eqref{eq:counterAssumption2}, next note that  
$(\ol{P},\ol{Q}')$ satisfies
\begin{align}
    & \ol{P} = 
    \l[
    \begin{array}{cc}
        U S_P & 0_{d_y\times 1} 
    \end{array}
    \r] \in \R^{d_y\times r} \nonumber\\
    &
     \ol{Q}' = 
     \l[
     \begin{array}{c}
     S_{Q'} V\\
     0_{1\times d_x}
     \end{array}
     \r] \in \R^{r\times m}
     ,\label{eq:contraductionAssumption}
\end{align}
where $U\in \R^{d_y\times (r-1)}$ and $V\in \R^{(r-1)\times m}$ correspond to those left and right singular vectors of $Y\proj_X$ that might be present in $\ol{P}\cdot \ol{Q}'$, see for example Lemma~5.1 (Item~5) in~\cite{hauser2018pca}. 

In~\eqref{eq:contraductionAssumption}, $S_P,S_{Q'} \in \R^{r\times r}$ are (not necessarily diagonal) matrices, and we note that $U$ and $V$ are column and row submatrices of $\wt{U}$ and $\wt{V}$, respectively. 

In view of~(\ref{eq:notRankDegRe},\ref{eq:svdYPX},\ref{eq:counterAssumption}), there exists a unique pair $(u,v)$ of left and right singular vectors of $Y\proj_X$ that is absent from~\eqref{eq:contraductionAssumption}, i.e.,
\begin{align}
    U^\top u = 0, \qquad V^\top v = 0.
    \label{eq:orthCounterAssumption}
\end{align}
To match the representation in~(\ref{eq:svdYPX}), note that $u$ above is a column-vector whereas $v$ is a row-vector. In particular,
\begin{align}
    & \wt{U} = 
    \l[
    \begin{array}{cc}
         U & u 
    \end{array} \r],
    \qquad 
    \wt{V} = 
    \l[
    \begin{array}{c}
         V  \\
         v 
    \end{array}
    \r], \nonumber\\
    & Y\proj_X = \wt{U}\wt{S} \wt{V} = U S V+usv,
    \qquad \text{(see \eqref{eq:svdYPX})}
    \label{eq:decompUVtildeProof}
\end{align}
where $S\in \R^{r-1}$ and $s\in \R$ collect the singular values corresponding to $(U,V)$ and $(u,v)$, respectively. 

To proceed, consider  inifinetsimally small scalars $\d_u$ and~$\d_v$.  Consider also an infinitesimally small perturbation $(\D_P,\D_{Q'})$ in $(\ol{P},\ol{Q}')$, specified as 
\begin{align}
    & \ol{P}+\D_P
     = \l[ 
    \begin{array}{cc}
        U S_P  & \d_u u
    \end{array}
      \r] \nonumber\\
          & \ol{Q}'+\D_{Q'}
     = \l[ 
    \begin{array}{c}
         S_{Q'} V  \\
         \d_v v
    \end{array}
      \r] .
      \label{eq:pertPrRankDeg}
\end{align}
It immediately follows that 
\begin{align}
    & (\ol{P}+\D_P)(\ol{Q}'+\D_{Q'}) \nonumber\\
    & =  \l[ 
    \begin{array}{cc}
        U S_P  & \d_u u
    \end{array}
      \r] 
      \cdot 
      \l[ 
    \begin{array}{c}
         S_{Q'} V  \\
         \d_v v
    \end{array}
      \r] \qquad \text{(see \eqref{eq:pertPrRankDeg})}  \nonumber\\
      & = U S_P S_{Q'} V + \d_u\d_v uv \nonumber\\
      & = \ol{P}\cdot \ol{Q}'+ \d_u\d_v uv.
      \qquad \text{(see \eqref{eq:contraductionAssumption})}
      \label{eq:prodPertRankDeg}
\end{align}
From~\eqref{eq:orthCounterAssumption}, it is evident that the perturbation in~\eqref{eq:prodPertRankDeg} is orthogonal to $\ol{P}\cdot \ol{Q}'$.

To continue, let us define the orthogonal projections $\proj_U=U U^\top$ and $\proj_u=uu^\top$, and define $\proj_V,\proj_v$ similarly. In particular, we can decompose $Y\proj_X$ as
\begin{align}
    Y\proj_X & = (\proj_U+\proj_u) (Y\proj_X) (\proj_V+\proj_v) \nonumber\\
    & = \proj_U (Y\proj_X) \proj_V + \proj_u (Y\proj_X) \proj_v,
    \label{eq:crossTermVanishY}
\end{align}
where the cross terms above vanish by properties of the SVD, see~(\ref{eq:svdYPX},\ref{eq:orthCounterAssumption}). Indeed, for example,
\begin{align}
     \proj_U (Y\proj_X) \proj_v & = UU^\top (\wt{U} \wt{S} \wt{V}) v^\top v 
     \qquad \text{(see \eqref{eq:svdYPX})}
     \nonumber\\
     & = UU^\top (U S V + u s v) v^\top v 
     \qquad \text{(see \eqref{eq:decompUVtildeProof})}
     \nonumber\\
     & = UU^\top  u s   v \qquad \text{(see \eqref{eq:orthCounterAssumption})} \nonumber\\
     & = 0, \qquad \text{(see \eqref{eq:orthCounterAssumption})}
\end{align}
where $S$ and $s$ collect the corresponding singular values for the singular vectors collected in $(U,V)$ and $(u,v)$, respectively. 
We will use the decomposition~\eqref{eq:crossTermVanishY} immediately below.

Under the  perturbation in~\eqref{eq:pertPrRankDeg}, the objective function of problem~\eqref{eq:pca} becomes
\begin{align}
    & \frac{1}{2} \|Y\proj_X - (\ol{P}+\D_P)(\ol{Q}'+\D_{Q'}) \|_F^2 \label{eq:decompPreRankDegenerate}\\
    & =\frac{1}{2} \|Y\proj_X - (\ol{P}+\d_u u)(\ol{Q}'+\d_v v) \|_F^2 
    \qquad \text{(see (\ref{eq:pertPrRankDeg}))}
    \nonumber\\
    & = \frac{1}{2}\| Y\proj_X - \ol{P}\cdot \ol{Q}' - \d_u\d_v uv\|_F^2 
    \qquad \text{(see \eqref{eq:prodPertRankDeg})}
    \nonumber\\
    & =\frac{1}{2}\| (\proj_U + \proj_u)  \nonumber\\
    & \qquad \cdot (Y\proj_X - \ol{P}\cdot \ol{Q}' - \d_u\d_v uv) (\proj_V + \proj_v) \|_F^2 \nonumber\\
    & = \frac{1}{2}\| (\proj_U (Y \proj_X) \proj_V - \ol{P}\cdot \ol{Q}') \nonumber\\
    & \qquad + (\proj_u (Y\proj_X) \proj_v -  \d_u\d_v uv) \|_F^2 
    \qquad \text{(see \eqref{eq:crossTermVanishY})}
    \nonumber\\
    & = \frac{1}{2}\| \proj_U (Y \proj_X) \proj_V - \ol{P}\cdot \ol{Q}'\|_F^2  \nonumber\\
    & \qquad + \frac{1}{2} \| \proj_u (Y\proj_X) \proj_v -  \d_u\d_v uv \|_F^2 \qquad \text{(see \eqref{eq:contraductionAssumption})} \nonumber\\
    & = \frac{1}{2}\| \proj_U (Y \proj_X) \proj_V - \ol{P}\cdot \ol{Q}'\|_F^2  \nonumber\\
    & \qquad + \frac{1}{2} | u^\top (Y\proj_X) v^\top -  \d_u\d_v  |^2.
    \label{eq:decreaseNotMin}
\end{align}    
It is now clear from~\eqref{eq:decreaseNotMin} that the perturbation in~\eqref{eq:pertPrRankDeg}  decreases the objective function of problem~\eqref{eq:pca} if we choose the signs of $\d_u$ and $\d_v$ carefully.

Indeed, we can upper bound the last line of~\eqref{eq:decreaseNotMin} as 
\begin{align} 
& \frac{1}{2} \|Y\proj_X - (\ol{P}+\d u)(\ol{Q}'+\d v) \|_F^2 \nonumber\\
        & = \frac{1}{2}\| \proj_U (Y \proj_X) \proj_V - \ol{P}\cdot \ol{Q}'\|_F^2  \nonumber\\
    & \qquad + \frac{1}{2} | u^\top (Y\proj_X) v^\top -  \d_u\d_v |^2 
    \qquad \text{(see \eqref{eq:decreaseNotMin})}
    \nonumber\\
    & < \frac{1}{2}\| \proj_U  (Y\proj_X) \proj_V - \ol{P}\cdot \ol{Q}'  \|_F^2 + \frac{1}{2} \| u^\top ( Y\proj_X) v \|_F^2 \nonumber\\
    & = \frac{1}{2}\| \proj_U  (Y\proj_X) \proj_V - \ol{P}\cdot \ol{Q}'  \|_F^2 + \frac{1}{2} \| \proj_u ( Y\proj_X) \proj_v \|_F^2 \nonumber\\
    & = \frac{1}{2}\|   Y\proj_X - \ol{P}\cdot \ol{Q}'  \|_F^2, \qquad \text{(see (\ref{eq:contraductionAssumption},\ref{eq:crossTermVanishY}))}
    \label{eq:decreaseNotMin2}
\end{align}
where we chose $\d_u\d_v$ above such that $\sign(\d_u\d_v)=\sign(u^\top(Y\proj_X)v^\top)$.

Note that~\eqref{eq:decreaseNotMin2}  contradicts the assumption that $(\ol{P},\ol{Q}')$ is an SOSP of problem~\eqref{eq:pca}, see Definition~\ref{defn:secondStN}. 
In fact, $(\ol{P},\ol{Q}')$ is a strict saddle point of problem~\eqref{eq:pca} because  $(\D_P,\D_{Q'})$ is a descent direction, see Definition~\ref{defn:strictSaddle}. 

Provided that $\rank(Y\proj_X)\ge r$, we conclude that any SOSP $(\ol{P},\ol{Q}')$ of problem~\eqref{eq:pca} satisfies 
\begin{align}
    \rank(\ol{W}) = r, \qquad \text{where}\,\, \ol{W} = \ol{P}\cdot \ol{Q}'.\label{eq:prodRankNonDeg}
\end{align}
We can in fact replace the conclusion in~\eqref{eq:prodRankNonDeg} with 
\begin{align}
    \rank(\ol{P}) = \rank(\ol{Q}') = \rank(\ol{W}) = r.\label{eq:indRankNonDeg}
\end{align}
Indeed, $\max(\rank(\ol{P}),\rank(\ol{Q}')) \le r$ because $\ol{P}$ has $r$ columns and $\ol{Q}'$ has $r$ rows, see~\eqref{eq:pca}. On the other hand, $\min(\ol{P},\ol{Q}')\ge r$ because of~\eqref{eq:prodRankNonDeg}. These two observations imply that $\rank(\ol{P})=\rank(\ol{Q}')=r$, as claimed in~\eqref{eq:indRankNonDeg}. 

Lastly, by EYM Theorem~\ref{thm:eym}, any SOSP of the PCA problem~\eqref{eq:pca} is also a global minimizer of problem~\eqref{eq:pca}. This completes the proof of Lemma~\ref{lem:noRankDegStPt}.

\section{Another Proof for Theorem~\ref{thm:wellBehavedMain}} \label{sec:anotherProofBah}

Here, we establish Theorem~\ref{thm:wellBehavedMain} with Proposition~32 in~\cite{bah2019learning} as the starting point. 
For completeness, let us first recall their result, adapted to our notation.
\begin{prop}
For every $r'\in [r]$, An SOSP $\ol{W}_{\NB}\in \M_{\NB,r'}$  of problem~\eqref{eq:mainN} is 
a global minimizer of $L_N$ restricted to the set $\M_{\NB,r'}$.
\end{prop}

Next, let us recall from~(\ref{eq:mainTwo},\ref{eq:pca}) that 
\begin{align}
    & \underset{W_{\NB}}{\min}\,\, L_N(W_{\NB}) \nonumber\\
    &= \underset{W_1,\cdots,W_N}{\min} \,\, \frac{1}{2}\| Y - W_N\cdots W_1 X\|_F^2 
    \qquad \text{(see (\ref{eq:mainN0},\ref{eq:mainN}))}
    \nonumber\\
    & = \underset{P,Q}{\min} \,\, \frac{1}{2}\|Y - PQX\|_F^2 
    \qquad \text{(see \eqref{eq:mainTwo})}
    \nonumber\\
    & = \frac{1}{2}\| Y\proj_{X^\perp}\|_F^2+ \min_{P,Q'}\,\, \frac{1}{2}\|Y\proj_X - PQ'\|_F^2,
    \label{eq:end-to-end}
\end{align}
where the last line above is from~\eqref{eq:pca}.

In view of~\eqref{eq:end-to-end}, a global minimizer $W_{\NB}^o = (W_1^o,\cdots, W_N^o)$ of problem~\eqref{eq:mainN} (the first program in~\eqref{eq:end-to-end}) corresponds to a global minimizer $(P^o,Q'^o)$ of problem~\eqref{eq:pca} (the last program in~\eqref{eq:end-to-end}) such that 
\begin{align}
    W_N^o\cdots W_1^o  =: W^o ,\qquad W^o X=  P^o Q'^o. \label{eq:globalO}
\end{align}
By assumption of Theorem~\ref{thm:wellBehavedMain}, it holds that $\rank(Y\proj_X)\ge r$. We can therefore invoke Lemma~\ref{lem:noRankDegStPt} in Appendix~\ref{sec:proofThmMain} to find that 
\begin{align}
    \rank(W_N^o\cdots W_1^o)
    & = \rank(W^o) \qquad \text{(see \eqref{eq:globalO})}  \nonumber\\
    & = r. \qquad \text{(see Lemma~\ref{lem:noRankDegStPt})}
    \label{eq:globalInMrN}
\end{align}
It is convenient to rewrite~\eqref{eq:globalInMrN} as
\begin{align}
(W_1^o,\cdots,W_N^o)\in \M_{\NB,r},
\label{eq:globalinMRNB}
\end{align}
where 
\begin{align}
     & \M_{\NB,r}:= \Big\{ W_{\NB}=(W_N,\cdots,W_1) \nonumber\\
     & \qquad \qquad \qquad : \rank(W_N\cdots W_1) = r \Big\} \subset \R^{d_{\NB}}. 
     \label{eq:defnMRNB}
\end{align}

On the other hand, with $k=r$, Proposition~32 in~\cite{bah2019learning} states that an  SOSP $\ol{W}_{\NB}\in \M_{\NB,r}$  of problem~\eqref{eq:mainN} is almost surely a global minimizer of $L_N$ restricted to the set $\M_{\NB,r}$. In view of~\eqref{eq:globalinMRNB}, we see that $\ol{W}_{\NB}$ is in fact a global minimizer of $L_N$ in $\R^{d_{\NB}}$. 
This completes our alternative proof for Theorem~\ref{thm:wellBehavedMain}.

\section{Theorem~35(a) in~\cite{bah2019learning} } \label{sec:revBahMain}

For completeness, here we recall Theorem~35(a) in~\cite{bah2019learning}, adapted to our notation. 

\textbf{Theorem 35(a)}~\cite{bah2019learning}\textbf{.} \emph{Suppose that $X$ has full column-rank and that $\rank(YX^\dagger X)\ge r$. Then gradient flow~\eqref{eq:gradFlowN} converges to a global minimizer of $L_N$ restricted to the $\M_{\NB,r'}$ for some $r'\le r$ from any initialization outside of a subset with Lebesgue measure zero. }

The key drawback of Theorem~35(a) above is it cannot ensure the convergence of gradient flow~\eqref{eq:gradFlowN} to a global minimizer of $L_N$. For example, if $r'=0$ above, then gradient flow converges to the zero matrix, which is known to be a non-strict saddle point of problem~\eqref{eq:mainN} when $N\ge 2$.

\section{Proof of Lemma \ref{lem:gfCvgFOSP}}\label{sec:proofGfCvgFOSP}

On the one hand, note that the objective function $L_N(W_{\NB})$ of problem~\eqref{eq:mainN} is analytic in $W_{\NB}$. 

On the other hand, recall the assumption that $XX^\top$ is invertible. Then, regardless of initialization,  gradient flow~\eqref{eq:gradFlowN} is bounded, i.e., contained in a finite ball centered at the origin, see Step 1 in the proof of Theorem~11 in~\cite{bah2019learning}.

We can now invoke the Lojasiewicz' theorem, see for example  Theorem~10 in~\cite{bah2019learning} or~\cite{absil2005convergence,lojasiewicz1982trajectoires}, to conclude that gradient flow converges to an FOSP $\ol{W}_{\NB}$ of problem~\eqref{eq:mainN}, regardless of initialization.

Lastly, gradient flow~\eqref{eq:gradFlowN} avoids strict saddle points of $L_N$ for almost every initialization  $W_{\NB,0}$,   with respect to the Lebesgue measure in $\R^{d_{\NB}}$, see Theorem~4.1 in~\cite{lee2016gradient}.\footnote{Strictly speaking,Theorem~4.1 in~\cite{lee2016gradient} is for gradient descent with a sufficiently small step size. However, their claim also holds for the limit case of  gradient flow, as the step size of gradient descent goes to zero.}

We conclude that the limit point $\ol{W}_{\NB}$ of gradient flow~\eqref{eq:gradFlowN} is in fact an SOSP of problem~\eqref{eq:mainN}, for almost every initialization with respect to the Lebesgue measure in $\R^{d_{\NB}}$.
This completes the proof of Lemma~\ref{lem:gfCvgFOSP}.  

Part of this argument is identical to the one in Theorem~11 of~\cite{bah2019learning}.

\section{Proof of Lemma \ref{lem:rankCteFlow}}\label{sec:proofRankCteFlow}

In the SVD of $W(t)$ in~\eqref{eq:analyticSVDCVG}, we let $\{s_i(t)\}_{i=1}^{\min(d_y,d_x)}$ denote the singular values of $W(t)$ in no particular order, with the corresponding left and right singular vectors denoted by  $\{u_{i}(t),v_{i}(t)\}_i$, for every $t\ge 0$.  

On the one hand, the evolution of the singular values of $W(t)$ in~\eqref{eq:analyticSVDCVG} is described by Theorem~3 in~\cite{arora2019implicit} as
\begin{align}
    \dot{s}_{i}(t) & = - N s_{i}(t)^{2-\frac{2}{N}} \cdot u_{i}(t)^\top \nabla L_1(W(t)) v_{i}(t),
    \label{eq:singsSpecified}
\end{align}
for every $t\ge 0$, where $L_1$ was defined in~\eqref{eq:mainOne}. Moreover, since the network depth $N\ge 2$, $\{s_i(t)\}_i$ remain nonnegative for every $t\ge 0$.

On the other hand, the singular values of $W(t)$ are bounded, i.e., $\max_i\sup_t s_i(t)<\infty$. 

Indeed, if $XX^\top$ is invertible, then gradient flow~\eqref{eq:gradFlowN} is bounded regardless of initialization, see Step 1 in the proof of Theorem~11 in~\cite{bah2019learning}. Recall also that gradient flow~\eqref{eq:gradFlowN} and induced flow~\eqref{eq:flowW} are related through the map 
\begin{align}
& \R^{d_{\NB}} \rightarrow \R^{d_y\times d_x} \nonumber\\
& W_{\NB}=(W_1,\cdots,W_N)\rightarrow W=W_N\cdots W_1.
\end{align}
Consequently, induced flow~\eqref{eq:flowW} and a fortiori its singular values too are bounded. 

We finally apply Lemma~4 in~\cite{arora2019implicit} to~\eqref{eq:singsSpecified} and find that each $s_{i}(t)$ is either zero for all $t\ge 0$, or positive for all $t\ge 0$, provided that the network depth $N\ge 2$. 

In other words,  $\rank(W(t))$ is invariant with $t$, i.e.,  
\begin{align}
    \rank(W(t)) = \rank(W_0), \qquad \forall t\ge 0,
\end{align}
which completes the proof of Lemma~\ref{lem:rankCteFlow}.

\section{Proof of Lemma~\ref{lem:genericNet}}\label{sec:proofGenericNet}

It is easy to see that $\M_{\NB,r}$ is not a closed set for any integer~$r$. For example, one can construct a sequence of rank-$1$ matrices that converge to the zero matrix. 

For the second claim in Lemma~\ref{lem:genericNet}, the proof is by induction over the depth $N$  of the linear network.

For the base of induction, when $N=1$, note that  $W_{\NB}=W_1\in \R^{d_1\times d_0}$, see~\eqref{eq:linearNetMap}.
It now follows from~\eqref{eq:rankNet} that 
\begin{align}
    \min(d_0,d_1) = r. \quad \text{(see (\ref{eq:rankNet}))}
    \label{eq:baseInduction}
\end{align}
In turn, it follows from~\eqref{eq:baseInduction} that almost every $W_1$ is rank-$r$, with respect to the Lebesgue measure in $\R^{d_{\NB}}=\R^{d_y\times d_x}$.

For the step of induction, suppose that 
\begin{align}
    \rank(W) = r, \qquad \text{with } W=W_{N}\cdots W_1, \label{eq:step}
\end{align}
for almost every $W_{\NB}=(W_1,\cdots,W_N)$, with respect to the Lebesgue measure in $\R^{d_{\NB}}$. In particular, it follows from~\eqref{eq:step} that $\range(W)$ is almost surely an $r$-dimensional subspace in $\R^{d_N}$, i.e., 
\begin{align}
    \dim(\range(W)) = r, \qquad \text{almost surely}. 
    \label{eq:rangeStep}
\end{align}

Consider a generic matrix $W_{N+1}$, with respect to the Lebesgue measure in $\R^{d_{N+1}\times d_N}$. We distinguish two cases. 

In the first case, suppose that  $d_{N+1}\ge d_{N}$. Then,  $W_{N+1}\in \R^{d_{N+1}\times d_N}$ has a trivial null space almost surely. With $\null$ standing for null space of a matrix, it follows that 
\begin{align}
     \null( W_{N+1} W) = \null(W),
     \label{eq:easyCaseLebesgue}
\end{align}
and, consequently, 
\begin{align}
    & \rank( W_{N+1} W_N\cdots W_1)  \nonumber\\
    & = \rank(W_{N+1}W) \qquad \text{(see \eqref{eq:step})} \nonumber\\
    & = d_0 - \dim(\null(W_{N+1}W))\nonumber\\
    & = d_0 - \dim(\null(W)) 
    \qquad \text{(see \eqref{eq:easyCaseLebesgue})}
    \nonumber\\
    & = \rank(W) = r, \label{eq:fullRank1}
\end{align}
almost surely with respect to the Lebesgue measure in $\R^{d_{N+1}\times d_N}$. Above, the third and last lines  use the fundamental theorem of linear algebra. 

In the second case, suppose that  $d_{N+1}<d_N$. Then the null space of $W_{N+1}\in \R^{d_{N+1}\times d_N}$ is a generic $(d_N-d_{N+1})$-dimensional subspace  of $\R^{d_N}$. It follows that 
\begin{align}
    \dim(\null(W_{N+1})) = d_N - d_{N+1} \le d_N - r,
    \label{eq:bndNullStep}
\end{align}
where the inequality above holds by~\eqref{eq:rankNet}. 
Note that 
\begin{align}
    & \dim(\range(W))+ \dim(\null(W_{N+1}))  \nonumber\\
    & \le r + (d_N - r) \qquad \text{(see (\ref{eq:rangeStep},\ref{eq:bndNullStep}))}  \nonumber\\
    & = d_N.\label{eq:notHit}
\end{align}
Since $\null(W_{N+1})$ is a generic subspace in $\R^{d_N}$  that satisfies~\eqref{eq:notHit}, it  almost surely holds that
\begin{align}
    \range(W) \cap\null(W_{N+1}) = \{0\}.\label{eq:noIntersectionStep}
\end{align}
Consequently,
\begin{align} 
    \null(W_{N+1}W) = \null(W), \qquad\text{(see \eqref{eq:noIntersectionStep})}
\end{align}
and it follows identically to~\eqref{eq:fullRank1} that 
\begin{align}
\rank(W_{N+1}\cdots W_1) = r,
    \label{eq:fullRank2}
\end{align}
almost surely with respect to the Lebesgue measure in~$\R^{d_{\NB+1}}$.  

We conclude from~(\ref{eq:fullRank1},\ref{eq:fullRank2}) that the  induction is complete, and this in turn completes the proof of the second and final claim in Lemma~\ref{lem:genericNet}.

\section{Proof of Lemma \ref{lem:remainNeighSimpleNBDomain}}\label{sec:proofRemainNeighSimpleNBDomain}

Recall that  the initialization of gradient flow~\eqref{eq:gradFlowN} is balanced by Assumption~\ref{assump:key2} and consider induced flow~\eqref{eq:flowW}. Let us define the set 
\begin{align}
    & \neigh_{\a}(Z_1)  := \Big\{ 
    W \overset{\svd}{=} u_W \cdot s_W \cdot v_W^\top : \nonumber\\
    &  \qquad s_W > (\a - \g_Z) s_Z, \nonumber\\ 
    & \qquad  u_W^\top Z_1 v_W > \a s_Z
    \Big\} \subset \R^{d_y\times d_x},
    \label{eq:neighWRe}
\end{align}
for $\a\in [\g_Z,1)$. 
Once initialized in $\neigh_{\a}(Z_1)$, induced flow remains there, as detailed in the next technical lemma.
\begin{lem}\label{lem:remainNeighSimple}
For induced flow~\eqref{eq:flowW} and $\a\in [\g_Z,1)$, 
$W_0\in \neigh_{\a}(Z_1)$ implies that $W_t\in \neigh_{\a}(Z_1)$ for all $t\ge 0$. That is, 
\begin{align}
    W_0 \in \neigh_\a(Z_1) \Longrightarrow W_t\in \neigh_\a(Z_1), \qquad \forall t\ge 0.
\end{align}
Above, Assumption~\emph{\ref{assump:key2}} and the notation therein are in force. 

\end{lem}
Before proving Lemma~\ref{lem:remainNeighSimple} in the next appendix, we show how it helps us prove Lemma~\ref{lem:remainNeighSimpleNBDomain}.

Indeed,  from Lemma~\ref{lem:remainNeighSimple} and the balanced initialization of gradient flow~\eqref{eq:gradFlowN}, it follows that 
\begin{align}
    W_{\NB,0}\in \neigh_{\NB,\a} \Longrightarrow W_{\NB,t} \in \neigh_{\NB,\a}, \qquad \forall t\ge 0,
    \label{eq:remainNBtemp}
\end{align}
under Assumption~\ref{assump:key2} and for $\a\in [\g_Z,1)$, where we used the definition of $\neigh_{\NB,\a}$ in~\eqref{eq:neighWNB}. 

Recall that the limit point $\ol{W}_{\NB}$ of gradient flow~\eqref{eq:gradFlowN} exists by Lemma~\ref{lem:gfCvgFOSP}, since $X$ has full-column rank by Assumption~\ref{assump:key2}. 

A byproduct of~\eqref{eq:remainNBtemp} about the limit point $\ol{W}_{\NB}$ of gradient flow~\eqref{eq:gradFlowN} is that 
\begin{align}
    & W_{\NB,0}\in \neigh_{\NB,\a} \Longrightarrow \nonumber\\
    & \qquad \ol{W}_{\NB} \in \text{closure}(\neigh_{\NB,\a}) \subset  \M_{\NB,1},
    \label{eq:remainsManifoldRdNReal}
\end{align}
where the set inclusion above holds true provided that 
$\a\in (\g_Z,1)$, see~(\ref{eq:preimage},\ref{eq:neighWNB},\ref{eq:remainNBtemp}). In words,~\eqref{eq:remainsManifoldRdNReal}  indicates that gradient flow does \emph{not} converge to the zero matrix. 

This completes the proof of Lemma~\ref{lem:remainNeighSimpleNBDomain}.

\subsection{Proof of Lemma \ref{lem:remainNeighSimple}}\label{sec:proofRemainNeighSimple}

The proof relies on the following technical lemma, which roughly-speaking states that the (rank-$1$) induced flow~\eqref{eq:flowW} always points in a similar direction as the (rank-$1$) target matrix $Z_1$.

\begin{lem}\label{lem:technicalNeigh}
Under Assumption~\emph{\ref{assump:key2}} and for $\a\in [\g_Z,1)$, $u_0^\top Z_1 v_0 > \a s_{Z}$ implies that $u_t^\top Z_1 v_t > \a s_{Z}$ for every $t\ge 0$. That is,
\begin{align}
    u_0^\top Z_1 v_0 >\a s_Z \Longrightarrow u_t^\top Z_1 v_t > \a s_Z, \,\, \forall t\ge 0.
\end{align}
Above, $W_t \overset{\svd}{=} u_t s_t v_t^\top$ is the rank-$1$ induced flow in~\emph{(\ref{eq:flowW},\ref{eq:analyticTSVD})}, and $s_Z,\g_{Z}$ were defined in~\eqref{eq:simpler}. 
\end{lem}
Before proving Lemma~\ref{lem:technicalNeigh} in the next appendix, let us see how Lemma~\ref{lem:technicalNeigh} can be used to prove Lemma~\ref{lem:remainNeighSimple}. 

Let us fix $\a\in[\g_Z,1)$. 
If $W_0\overset{\svd}{=} u_0 s_0 v_0^\top \in \neigh_{\a}(Z_1)$, then $u_0^\top Z_1 v_0^\top > \a s_{Z}$ by definition of $\neigh_{\a}(Z_1)$ in~\eqref{eq:neighWRe}. Lemma~\ref{lem:technicalNeigh} then implies that 
\begin{align}
    u_t^\top Z_1 v_t > \a s_{Z}, \qquad \forall t\ge  0. \qquad 
    \label{eq:oneCndMetNew}
\end{align}
To prove Lemma~\ref{lem:remainNeighSimple}, by the way of contradiction, let $\tau>0$ be the first time that the induced flow~\eqref{eq:flowW} leaves the set $\neigh_{\a}(Z_1)$. It thus holds that   
\begin{align}
    s_{\tau} & = \a s_{Z} - s_{Z,2} \qquad \text{(see \eqref{eq:neighWRe})} \nonumber\\
    & < u_\tau^\top Z_1 v_\tau - s_{Z,2}, \qquad \text{(see \eqref{eq:oneCndMetNew})}
    \label{eq:pushPreNew}
\end{align}
where the first line above uses the continuity of $s_t$ as a function of time $t$. Indeed, we know $s_t$ to be an analytic function of $t$, see~\eqref{eq:analyticTSVD}.

On the other hand, let us recall the evolution of the nonzero singular value of  flow~\eqref{eq:flowW} from~\eqref{eq:singsSpecified}, which we repeat here for convenience:
\begin{align}
    \dot{s}_\tau & = - N s_\tau^{2-\frac{2}{N}} \cdot u_\tau^\top \nabla L_1(W_\tau) v_\tau.
    \qquad \text{(see \eqref{eq:singsSpecified})}
    \label{eq:singValsWOGrad}
\end{align}
Recalling the definition of $L_1$ from~\eqref{eq:mainOne} and the whitened data assumption in~\eqref{eq:whitenedData}, we simplify the above gradient  as
\begin{align}
    \nabla L_1(W_\tau) & = W_\tau XX^\top - Y X^\top
    \qquad \text{(see \eqref{eq:mainOne})} 
    \nonumber\\
    & = m (W_\tau - Z).
    \qquad \text{(see (\ref{eq:whitenedData},\ref{eq:defnZ}))}
    \label{eq:gradWhitenedZ}
\end{align}
Substituting~\eqref{eq:gradWhitenedZ} back into~\eqref{eq:singValsWOGrad} and using the thin SVD of $W_\tau$ in~\eqref{eq:analyticTSVD}, we write at
\begin{align}
   \dot{s}_\tau 
   & =  - mN s_\tau^{2-\frac{2}{N}} \cdot u_\tau^\top (W_\tau - Z) v_\tau
   \qquad \text{(see (\ref{eq:singValsWOGrad},\ref{eq:gradWhitenedZ}))}
   \nonumber\\
   & =  - mN s_\tau^{2-\frac{2}{N}} \cdot ( s_\tau - u_\tau^\top   Z v_\tau) 
   \qquad \text{(see \eqref{eq:analyticTSVD})}
   \nonumber\\
   & > - mN s_\tau^{2-\frac{2}{N}} ( u_\tau^\top Z_1 v_\tau - s_{Z,2}  - u_\tau^\top Z v_\tau) 
   \,\, \text{(see \eqref{eq:pushPreNew})}
   \nonumber\\
   & = - mN s_\tau^{2-\frac{2}{N}} \big( u_\tau^\top Z_1 v_\tau - s_{Z,2} \nonumber\\
   & \qquad\qquad  - u_\tau^\top Z_1 v_\tau - u_\tau^\top Z_{1^+} v_\tau \big)
   \qquad \text{(see (\ref{eq:decomZ}))} \nonumber\\
   & =  - mN s_\tau^{2-\frac{2}{N}} ( - s_{Z,2} - u_\tau^\top Z_{1^+} v_\tau ) \nonumber\\
    &\ge 0, 
    \label{eq:pushedUp}
\end{align}
which pushes the singular value up and thus pushes the induced flow back into $\neigh_{\a}(Z_1)$. That is, the induced flow cannot escape from $\neigh_{\a}(Z_1)$. 

In the last line of~\eqref{eq:pushedUp}, we used the fact that $s_{Z,2}$ is the second largest singular value of $Z$ and hence the largest singular value of the residual matrix $Z_{1^+}$, see~\eqref{eq:specZRecallApp}. In the same line, we also used the fact that $u_\tau,v_\tau$ are unit-length vectors by construction, so that $u_\tau^\top Z_{1^+} v_\tau \ge - s_{Z,2}$.  
This completes the proof of Lemma~\ref{lem:remainNeighSimple}.

\subsection{Proof of Lemma \ref{lem:technicalNeigh}}
From~\eqref{eq:analyticTSVD}, recall the thin SVD of  induced flow~\eqref{eq:flowW}, i.e.,  
\begin{align}
W_t\overset{\svd}{=} u_t s_t v_t^\top, \qquad \forall t\ge 0,
\label{eq:thinSVDrecall}
\end{align}
where 
\begin{align}
\|u_t\|_2^2=\|v_t\|_2^2=1,
\label{eq:normCte}
\end{align}
and the only nonzero singular value is  $s_t> 0$. 

By taking the derivative of the  identities in~\eqref{eq:normCte} with respect to $t$, we find that 
\begin{align}
    & u_t^\top \dot{u}_t = 0, \nonumber\\
    & v_t^\top \dot{v}_t = 0, \qquad \forall t\ge0.
    \label{eq:alongDers}
\end{align}
By taking derivative of both sides of the thin SVD~\eqref{eq:thinSVDrecall}, we also find that
\begin{align}
    \dot{W}_t = \dot{u}_t s_t v_t^\top + u_t \dot{s}_t v_t^\top + u_t s_t \dot{v}_t^\top , \qquad \forall t\ge 0. \label{eq:derThinSVD}
\end{align}
Let $U_t\in \R^{d_y\times (d_y-1)}$ with orthonormal columns be orthogonal to  $u_t$.  By multiplying both sides of~\eqref{eq:derThinSVD} by $U_t^{\top}$, we find that 
\begin{align}
    U_t^{\top} \dot{W}_t & = U_t^{\top} \dot{u}_t s_t v_t^\top,  \qquad \forall t\ge 0,
\end{align}
which after rearranging yields that 
\begin{align}
    U_t^{\top} \dot{u}_t = s_t^{-1}U_t^{\top} \dot{W}_t v_t, \qquad \forall t\ge 0.
    \label{eq:deruDotRaw}
\end{align}
Combining~(\ref{eq:alongDers},\ref{eq:deruDotRaw}) yields that 
\begin{align}
    \dot{u}_t = s_t^{-1} \proj_{U_t} \dot{W}_t v_t, \qquad \forall t\ge 0,
    \label{eq:dotuFullRaw}
\end{align}
where $\proj_{U_t}=U_tU_t^\top$ is the orthogonal projection onto the subspace orthogonal to $u_t$.

Similarly, let  $V_t\in \R^{d_x\times (d_x-1)}$ with orthonormal columns be orthogonal to $v_t$. 
As before, by multiplying both sides of~\eqref{eq:derThinSVD} by $V_t$, we find that 
\begin{align}
\dot{W}_t V_t & = u_t s_t \dot{v}_t ^\top V_t, \qquad \forall t\ge 0,
\end{align}
which after rearranging yields
\begin{align}
    V_t^{\top} \dot{v}_t =  s_t^{-1} V_t^{\top} \dot{W}_t^\top u_t, \qquad \forall t\ge 0. 
    \label{eq:dervDotRaw}
\end{align}
Then, combining~(\ref{eq:alongDers},\ref{eq:dervDotRaw}) leads us to 
\begin{align}
    \dot{v}_t & = s_t^{-1} \proj_{V_t} \dot{W}_t^\top u_t, \qquad \forall t\ge 0,
    \label{eq:dotvFullRaw}
\end{align}
where $\proj_{V_t}=V_tV_t^\top$.

Both expressions~(\ref{eq:dotuFullRaw},\ref{eq:dotvFullRaw}) involve $\dot{W}_t$. Under the assumption of whitened data in~\eqref{eq:whitenedData}, we express $\dot{W}_t$ as
\begin{align}
    \dot{W}_t
    & = - \A_{W_t} ( W_t XX^\top - YX^\top ) 
    \qquad \text{(see \eqref{eq:flowW})} \nonumber\\
    &= - m \A_{W_t}(W_t -Z) \qquad \text{(see (\ref{eq:whitenedData},\ref{eq:defnZ}))} \nonumber\\
    & = - m N s_t^{1-\frac{2}{N}} (s_t - u_t^\top Z v_t) W_t 
    \qquad \text{(see (\ref{eq:defnAZ}))}
    \nonumber\\
    & \qquad +ms_t^{2-\frac{2}{N}}  \proj_{u_t} Z \proj_{V_t}  \nonumber\\
    & \qquad + ms_t^{2-\frac{2}{N}}  \proj_{U_t} Z \proj_{v_t}, \qquad \forall t\ge 0,
    \label{eq:expWDot}
\end{align}
where $\proj_{u_t}= u_tu_t^\top$ and $\proj_{v_t}=v_tv_t^\top$. The last identity above invokes the first part of Lemma~\ref{lem:innProdBnds}, which collects some basic properties of the operator $\A_W$. 

Substituting $\dot{W}_t$ back into~(\ref{eq:dotuFullRaw},\ref{eq:dotvFullRaw}), we reach
\begin{align}
    \dot{u}_t & = ms_t^{1-\frac{2}{N}} \proj_{U_t} Z v_t,  
    \qquad \text{(see \eqref{eq:expWDot})}
    \nonumber\\
    \dot{v}_t & = ms_t^{1-\frac{2}{N}} \proj_{V_t} Z^\top u_t, \qquad \forall t\ge 0.
    \label{eq:uvDotIndFinal}
\end{align}
It immediately follows from the first identity in~\eqref{eq:uvDotIndFinal} that  
\begin{align}
    u_Z^\top \dot{u}_t & = ms_t^{1-\frac{2}{N}} u_Z^\top \proj_{U_t} Z v_t
    \qquad \text{(see \eqref{eq:uvDotIndFinal})}
    \nonumber\\
    & = ms_t^{1-\frac{2}{N}} u_Z^\top \proj_{U_t} (u_Z s_Z v_Z^\top + Z_{1^+} ) v_t \,\, \text{(see \eqref{eq:decomZ})}
    \nonumber\\
    & = ms_t^{1-\frac{2}{N}} s_Z u_Z^\top \proj_{U_t} u_Z \cdot v_Z^\top v_t \nonumber\\
    & \qquad +m s_t^{1-\frac{2}{N}} u_Z^\top \proj_{U_t} Z_{1^+} v_t \nonumber\\
    & = ms_t^{1-\frac{2}{N}} s_Z \| \proj_{U_t} u_Z\|_2^2 \cdot v_Z^\top v_t \nonumber\\
    & \qquad + m s_t^{1-\frac{2}{N}} u_Z^\top \proj_{U_t} Z_{1^+} v_t, \quad \forall t\ge 0.
    \label{eq:leg1Remain}
\end{align}
To bound the last term above, note that 
\begin{align}
    & |u_Z^\top \proj_{U_t} Z_{1^+} v_t| \nonumber\\
    & \le  \| \proj_{U_t} u_Z \|_2 \cdot \| Z_{1^+} v_t\|_2 
    \qquad \text{(Cauchy-Schawrz ineq.)}
    \nonumber\\
    & \le \| \proj_{U_t} u_Z \|_2 \cdot s_{Z,2} \| \proj_{V_Z} v_t\|_2,
    \,\,\, \text{(see (\ref{eq:decomZ},\ref{eq:specZRecallApp}))}
        \label{eq:leg2Remain}
\end{align}
where $s_{Z,2}$ is the second largest singular value of $Z$ and thus the largest singular value of the residual matrix $Z_{1^+}$. Above,  $V_Z\in \R^{d_x\times (d_x-1)}$ with orthonormal columns  is orthogonal to $v_Z$, see~(\ref{eq:decomZ},\ref{eq:specZRecallApp}).

Similarly, it follows from the second identity in~\eqref{eq:uvDotIndFinal} that 
\begin{align}
    v_Z^\top \dot{v}_t & =m s_t^{1-\frac{2}{N}} v_Z^\top \proj_{V_t} Z^\top u_t 
    \qquad \text{(see \eqref{eq:uvDotIndFinal})}
    \nonumber\\
    & = ms_t^{1-\frac{1}{N}} v_Z^\top \proj_{V_t}( v_Z s_Z u_Z^\top  +Z_{1^+}^\top) u_t 
    \,\, \text{(see \eqref{eq:decomZ})}
    \nonumber\\
    & = ms_t^{1-\frac{2}{N}}s_Z \| \proj_{V_t} v_Z\|_2^2 \cdot u_Z^\top u_t \nonumber\\
    & \qquad + ms_t^{1-\frac{2}{N}}  v_Z^\top \proj_{V_t} Z_{1^+}^\top u_t, \quad \forall t\ge 0.
    \label{eq:leg3Remain}
\end{align}
To bound the last term above, we write that 
\begin{align}
    & |v_Z^\top \proj_{V_t} Z_{1^+}^\top u_t| \nonumber\\
    & \le \| \proj_{V_t} v_Z\|_2 \cdot \| Z_{1^+}^\top u_t\|_2 \nonumber\\
    & \le \| \proj_{V_t} v_Z\|_2 \cdot s_{Z,2} \| \proj_{U_Z} u_t\|_2,
    \label{eq:leg4Remain}
\end{align}
where $U_Z\in \R^{d_y\times (d_y-1)}$ with orthonormal columns is orthogonal to $u_Z$, see~\eqref{eq:decomZ}.

All these calculations in~(\ref{eq:leg1Remain}-\ref{eq:leg4Remain}) allow us to write that 
\begin{align}
    & \frac{\der (u_t^\top Z_1 v_t)}{\der t} \nonumber\\
    & = s_Z \frac{\der (u_t^\top u_Z v_Z^\top v_t)}{\der t} 
    \qquad \text{(see \eqref{eq:decomZ})}
    \nonumber\\
    & = s_Z (u_Z^\top \dot{u}_t) ( v_Z^\top v_t ) \nonumber\\
    &\qquad  + s_Z (u_Z^\top u_t) (v_Z^\top \dot{v}_t) 
     \,\, \text{(product rule)}
     \nonumber\\
    & = ms_Z^2 s_t^{1-\frac{2}{N}} \|\proj_{U_t} u_Z\|_2^2 \cdot (v_Z^\top v_t)^2\nonumber\\
    & \qquad + ms_Z^2 s_t^{1-\frac{2}{N}} \|\proj_{V_t} v_Z\|_2^2 \cdot (u_Z^\top u_t)^2 \nonumber\\
    & \qquad+ R_t,  \qquad \forall t\ge 0,
    \label{eq:prodDerComplete}
\end{align}
where the residual $R_t$ satisfies 
\begin{align}
    |R_t| & \le m s_Z s_{Z,2} s_t^{1-\frac{2}{N}} \|\proj_{U_t}u_Z\|_2 \|\proj_{V_Z}v_t\|_2    |v_Z^\top v_t| \nonumber\\
    & + m s_Z s_{Z,2} s_t^{1-\frac{2}{N}} \|\proj_{V_t}v_Z\|_2 \|\proj_{U_Z} u_t\|_2   |u_Z^\top u_t|.
    \label{eq:residialOriginal}
\end{align}
Let us set 
\begin{align}
    a_t := u_Z^\top u_t, \qquad b_t := v_Z^\top v_t,
    \label{eq:defnAtBt}
\end{align}
for short. Then we can rewrite~(\ref{eq:prodDerComplete},\ref{eq:residialOriginal}) as 
\begin{align}
     \frac{\der (u_t^\top Z_1 v_t)}{\der t} 
    & = ms_Z^2 s_t^{1-\frac{2}{N}} \l( (1- a_t^2)b_t^2 + a_t^2 (1-b_t^2)  \r) \nonumber\\
    & \qquad +R_t,\qquad \forall t\ge 0, 
    \label{eq:prodDerCompleteRewritten}
\end{align}
where 
\begin{align}
    |R_t| & \le ms_Z s_{Z,2} s_t^{1-\frac{2}{N}} \sqrt{1-a_t^2}\sqrt{1-b_t^2} \l(|a_t| + |b_t|    \r) \nonumber\\
    & \le  2m s_Z s_{Z,2} s_t^{1-\frac{2}{N}} \sqrt{(1-a_t^2)(1-b_t^2)} \nonumber\\
    & \le 2m s_Z s_{Z,2} s_t^{1-\frac{2}{N}} (1-a_tb_t),
    \label{eq:RinABnotation}
\end{align}
and the second  above uses the fact that 
\begin{align}
    |a_t| = |u_Z^\top u_t| \le \|u_Z\|_2\cdot \|u_t\|_2 \le 1,
    \label{eq:bndOnAtBt}
\end{align}
for every $t\ge 0$, and similarly $|b_t|\le 1$. The third line in~\eqref{eq:RinABnotation} uses the inequality
\begin{align}
    \sqrt{(1-a_t^2)(1-b_t^2)} & 
     = \sqrt{1-a_t^2-b_t^2 + a_t^2 b_t^2} \nonumber\\
     & \le \sqrt{1- 2a_tb_t + a_t^2 b_t^2} \nonumber\\
     & = 1-a_tb_t,
\end{align}
where the last line above again uses~\eqref{eq:bndOnAtBt}.  

The residual $R_t$ is small when the spectral gap of $Z$ is large. 
Indeed, note that 
\begin{align}
    |R_t| & \le 2m s_Z s_{Z,2} s_t^{1-\frac{2}{N}} (1-a_tb_t)
    \quad \text{(see \eqref{eq:RinABnotation})}\nonumber\\
    & < 2m s_Z^2 s_t^{1-\frac{2}{N}} a_tb_t (1-a_tb_t),
    \label{eq:simplifyR1}
\end{align}
where the last line above  holds provided that 
\begin{align}
    a_tb_t & > \frac{s_{Z,2}}{s_Z} 
     = \g_Z. \qquad \text{(see \eqref{eq:simpler})}
    \label{eq:specSmallAB}
\end{align}
We continue and bound the last line of~\eqref{eq:simplifyR1} as
\begin{align}
    |R_t| & < 2m s_Z^2 s_t^{1-\frac{2}{N}} a_tb_t (1-a_tb_t)
    \qquad \text{(see \eqref{eq:simplifyR1})}
    \nonumber\\
    & = 2m s_Z^2 s_t^{1-\frac{2}{N}}  ( a_tb_t-a_t^2 b_t^2) \nonumber\\
    & \le  2m s_Z^2 s_t^{1-\frac{2}{N}}  \l( \frac{a_t^2+b_t^2}{2} -a_t^2 b_t^2\r) \nonumber\\
    & =ms_Z^2 s_t^{1-\frac{2}{N}} (a_t^2+b_t^2 - 2 a_t^2 b_t^2 ) \nonumber\\
    & = ms_Z^2 s_t^{1-\frac{2}{N}} ( (1-a_t^2)b_t^2 + a_t^2 (1-b_t^2)). 
\end{align}
By comparing the above bound on the residual $R_t$ with~\eqref{eq:prodDerCompleteRewritten}, for a fixed time $t$, we conclude that
\begin{align}
    \frac{\der (u_t^\top Z_1 v_t )}{\der t} > 0,
\end{align}
provided that 
\begin{align}
    u_t^\top Z_1 v_t &  = s_Z \cdot u_t^\top u_Z v_Z^\top v_t \qquad \text{(see \eqref{eq:decomZ})} \nonumber\\
    & = s_Z \cdot a_t b_t \qquad \text{(see \eqref{eq:defnAtBt})}  \nonumber\\
    & > s_Z \g_Z.\qquad \text{(see \eqref{eq:specSmallAB})}
    \label{eq:abLarge}
\end{align}
For $\a\in [\g_Z,1)$, it immediately follows from~\eqref{eq:abLarge} that 
\begin{align}
    u_{0}^\top Z_1 v_0 > \a s_{Z} \Longrightarrow u_{t}^\top Z_1 v_t > \a s_{Z},
\end{align}
for every $t\ge 0$, which completes the proof of Lemma~\ref{lem:technicalNeigh}.

\section{Lazy Training}\label{sec:lazy}

For completeness, here we verify that Theorem~1 in~\cite{arora2018convergence} suffers from lazy training~\cite{chizat2019lazy}. To begin, let us recall Theorem 1 in~\cite{arora2018convergence}, adapted to our setting. 

\begin{thm}
    For $c>0$, suppose that gradient flow is initialized at $W_{\NB,0}=(W_{1,0},\cdots,W_{N,0})$ where $W_0 = W_{N,0}\cdots W_{1,0}$ satisfies $\|W_0 - Z\|_F \le \s_{\min}(Z) - c$. Here, $\s_{\min}(Z)$ stands for the smallest singular value of $Z$, defined in~\eqref{eq:defnZ}. Suppose also that $W_{\NB,0}$ is balanced, see Definition~\ref{defn:balanced}. Then, for a target accuracy $\epsilon>0$, it holds that 
    \begin{align}
        & l(t) =  \frac{1}{2}\| W(t) - Z\|_F^2 \le \epsilon, \nonumber\\
        &
         \qquad \forall t \ge \frac{1}{c^{2(1-\frac{1}{N})}} \log( l(0) / \epsilon),
    \end{align}
    where $W(t)= W_N(t)\cdots W_1(t)$. 
\end{thm}

For the sake of clarity, let us assume that $d_x=m$ and $X=\sqrt{m} I_{d_x}$, which satisfies the whitened requirement in~\eqref{eq:whitenedData} and~\cite{arora2018convergence}. Here, $I_{d_x}\in \R^{d_x\times d_x}$ is the identity matrix.

Recalling the loss function $L_N$ in~\eqref{eq:mainN} and the initialization $W_{\NB,0}\in \R^{d_{\NB}}$ of gradient flow~\eqref{eq:gradFlowN}, we write that 
\begin{align}
    & L_N(W_{\NB,0}) \nonumber\\
    & = \frac{1}{2}\| Y - W_{N,0}\cdots W_{1,0} X \|_F^2 \qquad \text{(see \eqref{eq:mainN})} \nonumber\\
    & = \frac{1}{2} \| Y - \sqrt{m} W_{N,0} \cdots W_{1,0} \|_F^2 \quad (X= \sqrt{m}I_{d_x})\nonumber\\
    & = \frac{m}{2} \l\| \frac{Y X^\top}{m} - W_{n,0}\cdots W_{1,0} \r\|_F^2 
    \quad (X= \sqrt{m}I_{d_x})
    \nonumber\\
    & = \frac{m}{2} \| Z - W_{N,0}\cdots W_{1,0}\|_F^2. \qquad  \text{(see \eqref{eq:defnZ})}
\end{align}
Definition~2 in~\cite{arora2018convergence} requires the last line above and, consequently, $L_N(W_{\NB,0})$ to be small. In turn, $L_N(W_{\NB,0})$ appears in Equation~(1) in~\cite{chizat2019lazy}.  Definition~2 in~\cite{arora2018convergence} thus requires the factor $\kappa$ in Equation~(1) in~\cite{chizat2019lazy} to be small, which is how the authors define the lazy training regime there.

\section{Proof of Lemma \ref{lem:switchFocusRate}} \label{sec:proofSwitchFocusRate}

From Theorem~\ref{cor:conjProved}, recall that gradient flow~\eqref{eq:gradFlowN} converges to a solution of~\eqref{eq:mainN} from almost every balanced initialization in the set $\neigh_{\NB,\a}$. That is,
\begin{align}
    & \lim_{t\rightarrow\infty} \frac{1}{2}\l\| Y - W_N(t)\cdots W_1(t)X\r\|_F^2 \nonumber\\
    & = \underset{W_1,\cdots,W_N}{\min} \frac{1}{2}\l\| Y - W_N\cdots W_1 X\r\|_F^2.
    \label{eq:cvgFlowObvious}
\end{align}

On the other hand, recall that gradient flow~\eqref{eq:gradFlowN} induces the flow~\eqref{eq:flowW} under the surjective map 
\begin{align}
& \R^{d_{\NB}} \rightarrow \M_{1,\cdots,r}\nonumber\\
& W_{\NB}=(W_1,\cdots,W_N)\rightarrow W=W_N\cdots W_1,
\label{eq:surjectiveMapProof}
\end{align}
where $\M_{1,\cdots,r}$ is the set of all $d_y\times d_x$ matrices  of rank at most $r$, see Appendix~\ref{sec:surjective} for the proof of the surjective property.

In view of~\eqref{eq:cvgFlowObvious}, induced flow~\eqref{eq:flowW} therefore satisfies 
\begin{align}
    & \lim_{t\rightarrow\infty} \frac{1}{2}\l\| Y - W(t)X\r\|_F^2 \nonumber\\
    & = \underset{\rank(W)\le r}{\min} \frac{1}{2}\l\| Y - W X\r\|_F^2,
    \label{eq:cvgInducedFlowObvious}
\end{align}
where $W(t)=W_N(t)\cdots W_1(t)$. 

Let $\proj_X=X^\dagger X$ and $\proj_{X^\perp}=I_m-\proj_X$ denote the orthogonal projections onto the row span of $X$ and its orthogonal complement, respectively. We can decompose $Y$ as
\begin{align}
    Y = Y\proj_X + Y \proj_{X}.
\end{align}
Using this decomposition, we can rewrite~\eqref{eq:cvgInducedFlowObvious} as
\begin{align}
        & \lim_{t\rightarrow\infty} \frac{1}{2}\l\| Y \proj_X - W(t)X\r\|_F^2 \nonumber\\
    & = \underset{\rank(W)\le r}{\min} \frac{1}{2}\l\| Y \proj_X - W X\r\|_F^2.
        \label{eq:cvgInducedFlowObvious2}
\end{align}
That is, in words, a linear network can only learn the component of $Y$ within the row span of $X$. 

Under Assumption~\ref{assump:key2}, the data matrix $X$ is whitened, so that $\proj_X = X^\dagger X = \frac{1}{m} X^\top X$, see~\eqref{eq:whitenedData}. We can therefore revise~\eqref{eq:cvgInducedFlowObvious2} as
\begin{align}
        & \lim_{t\rightarrow\infty} \frac{m}{2}\l\| Z - W(t)\r\|_F^2 \nonumber\\
     & = \underset{\rank(W)\le r}{\min} \frac{m}{2}\l\| Z - W \r\|_F^2,
        \label{eq:cvgInducedFlowObvious3}
\end{align}
where above we also used the definition of $Z$ in~\eqref{eq:defnZ}. 
To prove Lemma~\ref{lem:switchFocusRate}, we continue by setting $r=1$ in~\eqref{eq:cvgInducedFlowObvious3}. 

Recall also from Assumption~\ref{assump:key2} and specifically~\eqref{eq:simpler} that $Z$ has a nontrivial spectral gap, i.e., $s_{Z} >s_{Z,2}$. Therefore, $Z_1 =  u_Z s_Z v_Z^\top$ is the unique solution of the optimization problem in~\eqref{eq:cvgInducedFlowObvious3}, where the vectors $u_Z,v_Z$ are the corresponding leading left and right singular vectors of $Z_1$, see for example Section~1 in~\cite{golub1987generalization}. In view of this, it now follows from~\eqref{eq:cvgInducedFlowObvious3} with $r=1$ that 
\begin{align}
    \lim_{t\rightarrow \infty} \| Z_1 - W(t)\|_F = 0,
\end{align}
which completes the proof of Lemma~\ref{lem:switchFocusRate}.

\section{Derivation of (\ref{eq:lossToRHS})}\label{sec:derLossToRHS}

From Appendix~\ref{sec:proofEvolveLemma}, we will use the orthonormal bases $\wt{U}_t = [u_t,U_t]$ and $\wt{V}_t = [v_t,V_t]$. We decompose the loss function in these two bases as  
\begin{align}
    L_{1,1}(W_t) & = \frac{1}{2}\| W_t - Z_1\|_F^2 
    \qquad \text{(see \eqref{eq:newLossDefn})}
    \nonumber\\
    & = \frac{1}{2} \| \proj_{u_t}(W_t - Z_1)\proj_{v_t}\|_F^2 \nonumber\\
    & + \frac{1}{2} \| \proj_{u_t}(W_t - Z_1)\proj_{V_t}\|_F^2 \nonumber\\
    & + \frac{1}{2} \| \proj_{U_t}(W_t - Z_1)\proj_{v_t}\|_F^2 \nonumber\\
    & + \frac{1}{2} \| \proj_{U_t}(W_t - Z_1)\proj_{V_t}\|_F^2,
    \label{eq:newLossDecompose}
\end{align}
where $\proj_{u_t}=u_tu_t^\top$ is the orthogonal projection onto the span of $u_t$, and the remaining projection operators above are defined similarly.

Recalling the thin SVD $W_t = u_t s_t v_t^\top$ from~\eqref{eq:analyticTSVD} 
allows us to simplify~\eqref{eq:newLossDecompose} as
\begin{align}
    L_{1,1}(W_t) 
    & = \frac{1}{2} ( s_t - u_t^\top Z_1 v_t)^2 \nonumber\\
    & + \frac{1}{2} \| \proj_{u_t} Z_1\proj_{V_t}\|_F^2 + \frac{1}{2} \| \proj_{U_t} Z_1\proj_{v_t}\|_F^2  \nonumber\\
    & \qquad + \frac{1}{2} \| \proj_{U_t} Z_1\proj_{V_t}\|_F^2 \nonumber\\
    & = \frac{1}{2} ( s_t - u_t^\top Z_1 v_t)^2 \nonumber\\
    & \qquad + \frac{1}{2} \|Z_1\|_F^2 - \frac{1}{2} \| \proj_{u_t} Z_1\proj_{v_t} \|_F^2. 
    \label{eq:lossDecompose2}
\end{align}
Using the thin SVD $Z_1 = u_Z s_Z v_Z^\top$ from~\eqref{eq:decomZ} simplifies the last line above as
\begin{align}
    & \|Z_1\|_F^2 -  \| \proj_{u_t} Z_1\proj_{v_t} \|_F^2 \nonumber\\
    & = s_Z^2 - s_Z^2 (u_t^\top u_Z)^2 (v_Z^\top v_t)^2.
    \qquad \text{(see (\ref{eq:decomZ}))}
    \label{eq:decomposeZNorm}
\end{align}
Substituting the above identity back into~\eqref{eq:lossDecompose2} yields that
\begin{align}
   & L_{1,1}(W_t) \nonumber\\
   & = \frac{1}{2} ( s_t - u_t^\top Z_1 v_t)^2 \nonumber\\
   & \quad  + \frac{s_Z^2}{2} (1- (u_t^\top u_Z)^2 (v_Z^\top v_t)^2 ) 
   \,\,\text{(see (\ref{eq:lossDecompose2},\ref{eq:decomposeZNorm}))}
   \nonumber\\
   & = \frac{1}{2} ( s_t - u_t^\top Z_1 v_t)^2 \nonumber\\
   & \qquad + \frac{1}{2}(s_Z^2 - (u_t^\top Z_1 v_t)^2 )
      \qquad \text{(see \eqref{eq:decomZ})} \nonumber\\
     & = \frac{1}{2} ( s_t - u_t^\top Z_1 v_t)^2 \nonumber\\
   & \qquad + \frac{1}{2}(s_Z + u_t^\top Z_1 v_t ) (s_Z - u_t^\top Z_1 v_t )
       \nonumber\\
    &  \le   \frac{1}{2} ( s_t - u_t^\top Z_1 v_t)^2 \nonumber\\
   & \qquad  + s_Z (s_Z - u_t^\top Z_1 v_t ) ,
   \label{eq:uppBndL11Pr}
\end{align}
where the second identity and the inequality above use again the thin SVD of $Z_1$. The  inequality above also uses $u_t^\top Z_1 v_t \le s_Z$ twice, which holds true because $u_t,v_t$ are unit-length vectors by construction and $s_Z$ is the only nonzero singular of $Z_1$, see~\eqref{eq:decomZ}.  This completes the derivation of~\eqref{eq:lossToRHS}.

\section{Proof of Lemma~\ref{lem:evolveLemma}}\label{sec:proofEvolveLemma}

To begin, recall from~\eqref{eq:decomZ0} that $Z_1$ is the leading rank-$1$ component of $Z$, and let $Z_{1^+}= Z-Z_1$ denote the corresponding residual. We thus decompose $Z$ as 
\begin{align}
    Z & = Z_1 + Z_{1^+} \nonumber\\
    & \overset{\text{SVD}}{=} u_Z \cdot s_Z\cdot v_Z^\top  + U_Z S_Z V_Z^\top ,
    \label{eq:decomZ}
\end{align}
where $U_Z,V_Z$ contain the remaining left and right singular vectors of $Z$, and $S_Z$ contains the remaining singular values of of $Z$. In particular, let us repeat that
\begin{align}
    s_Z = \| Z\|=\|Z_1\|,\,\, s_{Z,2}=\|Z_{1^+}\|, \,\,\, \text{(see \eqref{eq:simpler})}
    \label{eq:specZRecallApp}
\end{align}
where $\|\cdot\|$ stands for spectral norm.

In this appendix, we compute the evolution of  loss function $L_{1,1}$ with time, which we recall from~\eqref{eq:evolveLossRaw} as
\begin{align}
    & \frac{\der L_{1,r}(W_t)}{\der t} \nonumber\\
    & = - m \langle W_t - Z_1, \A_{W_t}(W_t - Z) \rangle 
    \qquad \text{(see \eqref{eq:evolveLossRaw})} \nonumber\\
    & = - m \langle W_t - Z_1, \A_{W_t}(W_t - Z_1) \rangle  \nonumber\\
    & \qquad + m \langle W_t - Z_1, \A_{W_t}(Z_{1^+}) \rangle.
    \qquad \text{(see \eqref{eq:decomZ})}
    \label{eq:brkEvolveTwo}
\end{align}
To proceed, we will recall 
some basic properties of the operator $\A_W$, proved by algebraic manipulation of~\eqref{eq:defnAW} and included in  Appendix~\ref{sec:proofinnProdBnds} for completeness. The second identity below appears also in Lemma~5 of~\cite{bah2019learning}. 

\begin{lem}\label{lem:innProdBnds} For an arbitrary $W\in \R^{d_y\times d_x}$, let  
$$
W\overset{\text{SVD}}{=}\wt{U}\wt{S}\wt{V}^\top
$$
denote its  SVD, where $\wt{U},\wt{V}$ are orthonormal bases and $\wt{S}$ contains the singular values of $W$.  
Then, for an arbitrary $\D\in \R^{d_y\times d_x}$, it holds that 
\begin{align}
& \A_W(\D)  = \wt{U} \l(\sum_{j=1}^N (\wt{S}\wt{S}^\top)^{\frac{N-j}{N}} (\wt{U}^\top \D \wt{V})  (\wt{S}^\top \wt{S})^{\frac{j-1}{N}} \r) \wt{V}^\top,  \nonumber\\
&   \langle \D, \A_W(\D) \rangle\nonumber \\
& \qquad = \sum_{j=1}^N \l\| (\wt{S}\wt{S}^\top)^{\frac{N-j}{2N}} (\wt{U}^\top \D \wt{V}) (\wt{S}^\top \wt{S})^{\frac{j-1}{2N}}  \r\|_F^2. \label{eq:defnAZ} 
\end{align}
\end{lem}

For the first inner product in the last line of~\eqref{eq:brkEvolveTwo}, we invoke the second identity in~\eqref{eq:defnAZ} to write that 
\begin{align}
    & \langle W_t - Z_1, \A_{W_t}(W_t - Z_1) \rangle  \label{eq:firstTermEvolveApp}\\
    & = \sum_{j=1}^N \l\| (\wt{S}_t \wt{S}_t^\top)^{\frac{N-j}{2N}} \wt{U}_t^\top (W_t - Z_1) \wt{V}_t (\wt{S}_t^\top \wt{S}_t)^{\frac{j-1}{2N}} )   \r\|_F^2 
    \nonumber\\
    & =  \sum_{j=1}^N \l\| (\wt{S}_t \wt{S}_t^\top)^{\frac{N-j}{2N}} (\wt{S}_t - \wt{U}_t^\top   Z_1 \wt{V}_t) (\wt{S}_t^\top \wt{S}_t)^{\frac{j-1}{2N}} )   \r\|_F^2,\nonumber
\end{align}
where the second line above uses the SVD 
\begin{align}
W_t \overset{\text{SVD}}{=} \wt{U}_t \wt{S}_t \wt{V}_t^\top.
\qquad \text{(see \eqref{eq:analyticSVDCVG})}
\label{eq:fullSVDrecall}
\end{align}
We next  simplify the last line of~\eqref{eq:firstTermEvolveApp}.

In view of the thin SVD $W_t \overset{\svd}{=} u_t s_t v_t^\top$ in~\eqref{eq:analyticTSVD}, we let $U_{t}\in \R^{d_y\times (d_y-1)}$ and $V_t\in \R^{d_x\times (d_x-1)}$ be orthogonal complements for $u_t$ and $v_t$, respectively. This allows us to decompose $W_t$ as
\begin{align}
    W_t & \overset{\text{SVD}}{=} \wt{U}_t \wt{S}_t \wt{V}_t^\top 
    \qquad \text{(see \eqref{eq:fullSVDrecall})}
    \nonumber\\
    & = \l[ 
    \begin{array}{cc}
         u_t & U_t
    \end{array}  \r]
    \l[
    \begin{array}{cc}
         s_t &   \\
          & 0
    \end{array}
    \r]
    \l[
    \begin{array}{c}
         v_t^\top \\
         V_{t}^{\top}
    \end{array}
    \r],
    \label{eq:decompsWtSimple}
\end{align}
for every $t\ge 0$, where $0$ above is the $(d_y-1)\times (d_x-1)$ zero matrix. 
Using~\eqref{eq:decompsWtSimple}, we simplify~\eqref{eq:firstTermEvolveApp} to read 
\begin{align}
      & \langle W_t - Z_1, \A_{W_t}(W_t - Z_1) \rangle  \nonumber\\
      & = N s_t^{2-\frac{2}{N}} (s_t - u_t^\top Z_1 v_t)^2 \nonumber\\
      & \qquad + s_t^{2-\frac{2}{N}}\| u_t^\top Z_1 V_t\|_2^2\nonumber\\
     & \qquad + s_t^{2-\frac{2}{N}}\| U_t^{\top} Z_1 v_t\|_2^2.
     \label{eq:firstTermEvolveApp2}
\end{align}
The two norms above can be further simplified. Let us set
\begin{align}
    a_t := u_t^\top u_Z, \qquad b_t := v_t^\top v_Z, 
    \label{eq:abSecondTime}
\end{align}
for short.
Then we expand the first norm in~\eqref{eq:firstTermEvolveApp2} as
\begin{align}
    \|u_t^\top Z_1 V_t\|_2& = s_Z \| u_t^\top u_Z v_Z^\top V_t\|_2
    \qquad \text{(see \eqref{eq:decomZ})} \nonumber\\
    & = s_Z |u_t^\top u_Z| \cdot \| v_Z^\top V_t\|_2 \nonumber\\
    & = s_Z a_t \sqrt{1-b_t^2}, \qquad \text{(see (\ref{eq:abSecondTime}))}
    \label{eq:ThreeToTwo1}
\end{align}
where the last line above follows because $V_t$ spans the orthogonal complement of $v_t$. 
Likewise, the second norm in~\eqref{eq:firstTermEvolveApp2} is expanded as
\begin{align}
    \| U_t^\top Z_1 v_t\|_2 & = s_Z \|U_t^\top u_Z\|_2 \cdot |v_z^\top v_t| \nonumber\\
    & = s_Z b_t \sqrt{1-a_t^2} .
    \qquad \text{(see (\ref{eq:abSecondTime}))}
        \label{eq:ThreeToTwo2}
\end{align}
In particular, by combining~(\ref{eq:ThreeToTwo1},\ref{eq:ThreeToTwo2}), we find that 
\begin{align}
     & \|u_t^\top Z_1 V_t\|_2^2+\| U_t^\top Z_1 v_t\|_2^2 \nonumber\\
     &   =    s_Z^2 \l( a_t^2(1-b_t^2) +  (1-b_t^2)a_t^2  \r)
     \qquad \text{(see (\ref{eq:ThreeToTwo1},\ref{eq:ThreeToTwo2}))}
     \nonumber\\
      & =  s_Z^2 \l( a_t^2+b_t^2 - 2a_t^2 b_t^2 \r) \nonumber\\
     & \ge  s_Z^2 \l( 2a_tb_t - 2a_t^2 b_t^2 \r) \nonumber\\
      & = 2s_Z^2 a_tb_t( 1 - a_tb_t),
      \label{eq:ThreeToTwo2.1}
\end{align}
where the penultimate line above uses the inequality $a_t^2+b_t^2\ge 2a_tb_t$. 
Plugging~(\ref{eq:ThreeToTwo2.1}) back into~\eqref{eq:firstTermEvolveApp2}, we arrive at
\begin{align}
        & \langle W_t - Z_1, \A_{W_t}(W_t - Z_1) \rangle  \nonumber\\
      & \ge  N s_t^{2-\frac{2}{N}} (s_t - u_t^\top Z_1 v_t)^2 \nonumber\\
      & \qquad +2 s_t^{2-\frac{2}{N}} s_Z^2 
      a_tb_t( 1 - a_tb_t).
      \label{eq:firstTermEvolveApp2.5}
\end{align}
For the second inner product in the last line of~\eqref{eq:brkEvolveTwo}, we invoke the first identity in~\eqref{eq:defnAZ} to write that 
\begin{align}
    & \langle W_t - Z, \A_{W_t}(Z_{1^+}) \rangle \nonumber\\
    & = N s_t^{2-\frac{2}{N}} (s_t - u_t^\top Z_1 v_t) (u_t^\top Z_{1^+} v_t) \nonumber\\
    & \qquad - s_t^{2-\frac{2}{N}} \l\langle u_t^\top Z_1 V_t, u_t^\top Z_{1^+} V_t\r\rangle   \nonumber\\
    & \qquad - s_t^{2-\frac{2}{N}} \l\langle U_t^{\top} Z_1 v_t, U_t^{\top} Z_{1^+} v_t\r\rangle,  
\end{align}
and, consequently, 
\begin{align}
    & \l| \langle W_t - Z, \A_{W_t}(Z_{1^+}) \rangle\r| \nonumber\\
    & \le  N s_t^{2-\frac{2}{N}} |s_t - u_t^\top Z_1 v_t| \cdot |u_t^\top Z_{1^+} v_t| \nonumber\\
    & \qquad + s_t^{2-\frac{2}{N}} \|u_t^\top Z_1 V_t\|_2 \| u_t^\top Z_{1^+} V_t\|_2   \nonumber\\
    & \qquad + s_t^{2-\frac{2}{N}} \|U_t^{\top} Z_1 v_t\|_2 \|U_t^{\top} Z_{1^+} v_t\|_2, 
    \label{eq:secondTermEvolveApp}
\end{align}
where we twice used the Cauchy-Schwarz inequality above.

Recall the decomposition of $Z$ in~\eqref{eq:decomZ}. The three terms in~\eqref{eq:secondTermEvolveApp} that involve the residual matrix $Z_{1^+}$ can  be simplified as follows. For the first term, we write that 
\begin{align}
    |u_t^\top Z_{1^+} v_t| & = | u_t^\top U_Z S_Z V_Z^\top v_t|
    \qquad \text{(see \eqref{eq:decomZ})}\nonumber\\
    & \le  \| u_t^\top U_Z\|_2 \cdot \|S_Z\|\cdot  \| V_Z^\top v_t\|_2
    \nonumber\\
    & =s_{Z,2} \| u_t^\top U_Z\|_2 \cdot \| V_Z^\top v_t\|_2
    \nonumber\\
    & = s_{Z,2} \sqrt{1-a_t^2} \sqrt{1-b_t^2}
    \qquad \text{(see \eqref{eq:abSecondTime})} \nonumber\\
    & = s_{Z,2} \sqrt{1 - a_t^2 - b_t^2 +a_t^2 b_t^2} \nonumber\\
    & \le s_{Z,2} \sqrt{1 - 2 a_t  b_t +a_t^2 b_t^2} 
    \nonumber\\
    & = s_{Z,2} (1 -  a_t  b_t),
        \label{eq:ThreeToTwo2.5}
\end{align}
where $\|S_Z\|$ denotes the spectral norm of the matrix $S_Z$. The second line above uses the fact that $s_{Z,2}$ is the second largest singular value of $Z$ and thus the largest singular value of the residual matrix $Z_{1^+}$, see~(\ref{eq:decomZ},\ref{eq:specZRecallApp}).  The last line above uses the observation that $a_tb_t \le 1$ since $u_t,v_t,u_Z,v_Z$ all have unit norm by construction, see~\eqref{eq:abSecondTime}.

Likewise, another term in~\eqref{eq:secondTermEvolveApp} can be bounded as
\begin{align}
    \| u_t^\top Z_{1^+} V_t\|_2 & = \| u_t^\top U_Z S_Z V_Z^\top V_t\|_2 
    \qquad \text{(see \eqref{eq:decomZ})}
    \nonumber\\
    & \le \| u_t ^\top U_Z\|_2 \cdot \|S_Z\| \cdot \|V_Z^\top V_t\| \nonumber\\
    & = s_{Z,2}\| u_t ^\top U_Z\|_2 \cdot \|V_Z^\top V_t\| \nonumber\\
    & \le s_{Z,2}\| u_t ^\top U_Z\|_2 \cdot \|V_Z\| \cdot \| V_t\| \nonumber\\
    & \le s_{Z,2} \sqrt{1-a_t^2},
    \qquad \text{(see \eqref{eq:abSecondTime})}
    \label{eq:ThreeToTwo3}
\end{align}
where the last line above  uses the fact that both $V_Z$ and $V_t$ have orthonormal columns. 

For the last term involving $Z_{1^+}$ in~\eqref{eq:secondTermEvolveApp}, we similarly write that
\begin{align}
    \| U_t^\top Z_{1^+} v_t\|_2 & \le s_{Z,2} \| V_Z^\top v_t\|_2 = s_{Z,2} \sqrt{1-b_t^2}. 
    \label{eq:ThreeToTwo4}
\end{align}
Plugging back~(\ref{eq:ThreeToTwo1},\ref{eq:ThreeToTwo2},\ref{eq:ThreeToTwo2.5},\ref{eq:ThreeToTwo3},\ref{eq:ThreeToTwo4}) into~\eqref{eq:secondTermEvolveApp}, we reach 
\begin{align}
    & |\langle W_t - Z,\A_{W_t}(Z_{1^+}) \rangle| \nonumber\\
    & \le N s_t^{2-\frac{2}{N}} s_{Z,2} |s_t - u_t^\top Z_1 v_t| (1-a_tb_t) \nonumber\\
    & \qquad + s_{t}^{2-\frac{2}{N}} s_Z s_{Z,2} (a_t+b_t) \sqrt{1-a_t^2} \sqrt{1-b_t^2} \nonumber\\
    & \le N s_t^{2-\frac{2}{N}} s_{Z,2} |s_t - u_t^\top Z_1 v_t| (1-a_tb_t) \nonumber\\
    & \qquad + 2s_{t}^{2-\frac{2}{N}} s_Z s_{Z,2}  (1-a_tb_t), 
             \label{eq:secondTermEvolveApp2}
\end{align}
where the last line above follows from the chain of inequalities 
\begin{align}
    & (a_t+b_t) \sqrt{1-a_t^2}\sqrt{1-b_t^2} \nonumber\\
    & \le 2 \sqrt{1-a_t^2-b_t^2+a_t^2b_t^2} \qquad \text{(see \eqref{eq:abSecondTime})}\nonumber\\
    & \le 2 \sqrt{1-2a_tb_t+a_t^2b_t^2} \nonumber\\
    & = 2(1-a_tb_t).  \qquad \text{(see \eqref{eq:abSecondTime})}
\end{align}
Above, in the second and last lines, we used the fact that $a_t\le 1$ and $b_t\le 1$, see their definition in~\eqref{eq:abSecondTime}.

By combining~(\ref{eq:firstTermEvolveApp2.5},\ref{eq:secondTermEvolveApp2}), we can upper bound the evolution of the loss function as
\begin{align}
    & \frac{\der L_{1,1}(W_t)}{\der t} \nonumber\\
     & \le -m N s_t^{2-\frac{2}{N}} (s_t - u_t^\top Z_1 v_t)^2 \nonumber\\
      & \qquad -2 m s_t^{2-\frac{2}{N}} s_Z^2 a_tb_t(1-a_tb_t) \nonumber\\
      & \qquad + mNs_t^{2-\frac{2}{N}} s_{Z,2} |s_t - u_t^\top Z_1 v_t| (1-a_tb_t) \nonumber\\
      & \qquad + 2ms_t^{2-\frac{2}{N}} s_Z s_{Z,2} (1-a_tb_t) \nonumber\\
      & = -m N s_t^{2-\frac{2}{N}} (s_t - u_t^\top Z_1 v_t)^2 \nonumber\\
      & \qquad -2 m s_t^{2-\frac{2}{N}} ( u_t^\top Z_1 v_t) (s_Z-u_t^\top Z_1 v) \nonumber\\
      & \qquad + mNs_t^{2-\frac{2}{N}} \g_Z |s_t - u_t^\top Z_1 v_t| (s_Z -u_t^\top Z_1 v_t) \nonumber\\
      & \qquad + 2ms_t^{2-\frac{2}{N}}  s_{Z,2} (s_Z-u_t^\top Z_1 v_t), 
\end{align}
where the identity above uses the fact that  $s_Z a_t b_t = u_t^\top Z_1 v_t$, see~(\ref{eq:decomZ},\ref{eq:abSecondTime}). The inverse spectral gap $\g_Z=s_{Z,2}/s_Z$ was introduced in~(\ref{eq:simpler}). This completes the proof of Lemma~\ref{lem:evolveLemma}.

\subsection{Proof of Lemma~\ref{lem:innProdBnds}}\label{sec:proofinnProdBnds}

Let $W \overset{\text{SVD}}{=} \wt{U}\wt{S}\wt{V}^\top$ denote the  SVD of $W$, where $\wt{U},\wt{V}$ are orthonormal bases, and $\wt{S}$ contains the singular values of~$W$. 

Using the definition of $\A_W$, we write that 
\begin{align}
    \A_W(\D) & = \sum_{j=1}^N (W W^\top)^{\frac{N-j}{N}} \D (W^\top W)^{\frac{j-1}{N}}
    \quad \text{(see \eqref{eq:defnAW})}
    \nonumber\\
    & = \sum_{j=1}^N (\wt{U}\wt{S}\wt{S}^\top \wt{U}^\top)^{\frac{N-j}{N}} \D (\wt{V} \wt{S}^\top \wt{S} \wt{V}^\top)^{\frac{j-1}{N}}\nonumber\\
    & = \sum_{j=1}^N \wt{U} (\wt{S} \wt{S}^\top)^{\frac{N-j}{N}} \wt{U}^\top \D \wt{V} (\wt{S}^\top \wt{S})^{\frac{j-1}{N}} \wt{V}^\top \nonumber\\
    & =: \sum_{j=1}^N \A_{W,j}(\D),
    \label{eq:recallAW}
\end{align}
which proves the first claim in Lemma~\ref{lem:innProdBnds}.

For every $j\in \NB$, it also holds that 
\begin{align}
    & \langle \D,{\A}_{W,j}(\D) \rangle   \nonumber\\
    & = \langle \D, \wt{U} (\wt{S}\wt{S}^\top)^{\frac{N-j}{N}} \wt{U}^\top \D \wt{V} (\wt{S}^\top \wt{S})^{\frac{j-1}{N}} \wt{V}^\top \rangle 
    \quad \text{(see \eqref{eq:recallAW})}
    \nonumber\\
    & = \l\| (\wt{S}\wt{S}^\top)^{\frac{N-j}{2N}} \wt{U}^\top \D \wt{V} (\wt{S}^\top \wt{S})^{\frac{j-1}{2N}} \r\|_F^2,
    \label{eq:innFrob}
\end{align}
where the last line uses the fact that $\wt{U},\wt{V}$ are orthonormal bases. The proof of Lemma~\ref{lem:innProdBnds} is complete after summing up the above identity over $j$.

\section{Proof of Lemma \ref{lem:remainNeighSimpleTwo}}\label{sec:proofRemainNeighSimpleTwo}

The proof is  similar to that of Lemma~\ref{lem:remainNeighSimple}.

Let us fix $\a\in[\g_Z,1)$ and $\b>1$.  
If 
\begin{align}
W_0\overset{\svd}{=} u_0 s_0 v_0^\top \in \neigh_{\a,\b}(Z_1),
\end{align}
then $u_0^\top Z_1 v_0^\top > \a s_{Z}$ by definition of $\neigh_{\a,\b}(Z_1)$ in~\eqref{eq:neighWTwo}. Lemma~\ref{lem:technicalNeigh} then implies that 
\begin{align}
    u_t^\top Z_1 v_t > \a s_{Z}, \qquad \forall t\ge  0. \qquad 
    \label{eq:oneCndMet}
\end{align}
To prove Lemma~\ref{lem:remainNeighSimpleTwo}, by the way of contradiction, let $\tau>0$ be the first time that induced flow~\eqref{eq:flowW} leaves the set $\neigh_{\a,\b}(Z_1)$. It thus holds that   
\begin{align}
    s_{\tau} & = \a s_{Z} - s_{Z,2} ,
    \label{eq:pushPre}
\end{align}
or 
\begin{align}
    s_\tau = \b s_Z,
    \label{eq:pushOutPre}
\end{align}
where both of the identities above use the continuity of $s_t$ as a function of $t$. Indeed, we know $s_t$ to be an analytic function of time $t$, see~\eqref{eq:analyticTSVD}.

The case where~\eqref{eq:pushPre} happens is handled identically to the proof of Lemma~\ref{lem:remainNeighSimple}. We therefore focus on when the second case happens, i.e.,~\eqref{eq:pushOutPre}.

Recalling the second identity in~\eqref{eq:pushedUp}, we bound the evolution of the singular value of  induced flow~\eqref{eq:flowW} as
\begin{align}
   \dot{s}_\tau 
   & =  - mN s_\tau^{2-\frac{2}{N}} \cdot ( s_\tau - u_\tau^\top   Z v_\tau) 
   \qquad \text{(see \eqref{eq:pushedUp})}
   \nonumber\\
   & = - mN s_\tau^{2-\frac{2}{N}} ( \b s_{Z}  - u_\tau^\top Z v_\tau) 
   \qquad  \text{(see \eqref{eq:pushOutPre})}
   \nonumber\\
   & \le  - mN s_\tau^{2-\frac{2}{N}}( \b s_Z - s_Z)
   \qquad \text{(see (\ref{eq:decomZ}))} \nonumber\\
   & <0,    \label{eq:pushedUpTwo}
\end{align}
which pushes the singular value down and thus pushes the induced flow back into $\neigh_{\a,\b}(Z_1)$. That is, the induced flow cannot escape from $\neigh_{\a,\b}(Z_1)$. 

In the third line of~\eqref{eq:pushedUpTwo}, we used the fact that $s_Z$ is the leading singular value of $Z$, and $u_\tau,v_\tau$ are unit-norm vectors, see~(\ref{eq:decomZ},\ref{eq:thinSVDrecall}), thus $u_\tau^\top Z v_\tau \le s_Z$. The last line in~\eqref{eq:pushedUpTwo} holds because Lemma~\ref{lem:remainNeighSimpleTwo} assumes that $\b>1$.
This completes the proof of Lemma~\ref{lem:remainNeighSimpleTwo}.

\section{Proof of Theorem~\ref{thm:MainSimple}}\label{sec:proofMainSimple}

 Let us fix $\a\in [\g_Z,1)$ and $\b>1$. 
In view of Lemma~\ref{lem:remainNeighSimpleTwo}, we assume henceforth that induced flow~\eqref{eq:flowW} is initialized within $\neigh_{\a,\b}(Z_1)$ and thus remains there forever, i.e., 
\begin{align}
    W_t\in \neigh_{\a,\b}(Z_1),\qquad \forall t\ge 0.
    \label{eq:WstaysNew}
\end{align}

Using the definition of $\neigh_{\a,\b}(Z_1)$ in~(\ref{eq:neighWTwo}) and~(\ref{eq:WstaysNew}), we can update~\eqref{eq:evolveSimpleLemma} as
\begin{align}
        & \frac{\der L_{1,1}(W_t)}{\der t} \nonumber\\
     & \le  -  m N ((\a-\g_Z) s_Z)^{2-\frac{2}{N}} (s_t - u_t^\top Z_1 v_t)^2 
     \nonumber\\
      &  -2\a    m s_Z ((\a-\g_Z) s_Z)^{2-\frac{2}{N}}  (s_Z-u_t^\top Z_1 v_t)  \nonumber\\
      &  + m  N (\b s_Z)^{2-\frac{2}{N}} \g_Z |s_t - u_t^\top Z_1 v_t| (s_Z -u_t^\top Z_1 v_t) \nonumber\\
      & \qquad + 2 m (\b s_Z)^{2-\frac{2}{N}}  s_{Z,2} (s_Z-u_t^\top Z_1 v_t).
     \label{eq:evolveSimpleLemma2}
\end{align}
Recalling the upper bound on the loss function in~\eqref{eq:lossToRHS}, we can distinguish two regimes in the dynamics of~\eqref{eq:evolveSimpleLemma2}, depending on the dominant term on the right-hand side of~\eqref{eq:lossToRHS}, as detailed next.

\textbf{(Fast convergence)} When
\begin{align}
    \frac{1}{2}(s_t-u_t^\top Z_1 v_t)^2  \ge s_Z (s_Z - u_t^\top Z_1 v_t),\label{eq:fastCvg}
\end{align}
the loss function can be bounded as 
\begin{align}
    L_{1,1}(W_t) \le (s_t-u_t^\top Z_1 v_t)^2, \qquad \text{(see (\ref{eq:lossToRHS},\ref{eq:fastCvg}))}
    \label{eq:lossFastCvg}
\end{align}
and the evolution of loss in~\eqref{eq:evolveSimpleLemma2} thus simplifies to
\begin{align}
     \frac{\der L_{1,1}(W_t)}{\der t} 
      & \le -  m N s_Z^{2-\frac{2}{N}} \cdot    \label{eq:evolveFast}\\
      &  \l(  (\a-\g_Z)^{2-\frac{2}{N}}    - 2\g_Z\b^{2-\frac{2}{N}}  \r) \cdot  L_{1,1}(W_t) ,\nonumber
\end{align}
see Appendix~\ref{sec:derivationEvolveFastSlow} for the detailed derivation of~\eqref{eq:evolveFast}. 

\textbf{(Slow~convergence)} On the other hand, when 
\begin{align}
    \frac{1}{2}(s_t-u_t^\top Z_1 v_t)^2  \le  s_Z (s_Z - u_t^\top Z_1 v_t),\label{eq:slowCvg}
\end{align}
the loss function can be bounded as 
\begin{align}
L_{1,1}(W_t) \le 2 s_Z (s_Z - u_t^\top Z_1 v_t),
\,\, \text{(see (\ref{eq:lossToRHS},\ref{eq:slowCvg}))}
\label{eq:slowCvgLoss}
\end{align}
and the evolution of loss in~\eqref{eq:evolveSimpleLemma2} simplifies to
\begin{align}
      \frac{\der L_{1,1}(W_t)}{\der t} 
      & \le - m  s_Z^{2-\frac{2}{N}} \cdot \label{eq:evolveSlow} \\
      &  \l( \a (\a-\g_Z)^{2-\frac{2}{N}} - 2\g_Z N\b^{2-\frac{2}{N}}  \r) \cdot  L_{1,1}(W_t), \nonumber
\end{align}
see Appendix~\ref{sec:derivationEvolveFastSlow} again for the detailed derivation of~\eqref{eq:evolveSlow}.

In view of~(\ref{eq:fastCvg},\ref{eq:slowCvg}), the key transition between fast and slow convergence rates happens when 
\begin{align}
T_t& :=T_{1,t}-s_Z T_{2,t} \nonumber\\
& = \frac{1}{2}(s_t - u_t^\top Z_1 v_t)^2 - s_Z (s_Z - u_t^\top Z_1 v_t) 
\end{align}
changes sign. Above, we used the definition of $T_{1,t}$ and $T_{2,t}$ in~\eqref{eq:lossToRHS}. 

Instead of the first time such  a sign change happens, it is convenient to  consider the more conservative choice  of time $\tau \ge 0$ when
\begin{align}
    s_{\tau } \le \sqrt{6} s_Z
    \label{eq:conservative}
\end{align}
for the first time. 
Indeed, if~\eqref{eq:conservative} does not hold, then  $T_\tau > 0$ and thus the fast convergence is in force. This claim is verified in~Appendix~\ref{sec:derConservative} for  completeness. 

With the definition of $\tau$ at hand from~\eqref{eq:conservative}, we can combine~\eqref{eq:evolveFast} and~\eqref{eq:evolveSlow} to obtain that 
\begin{align}
    & \frac{\der L_{1,1}(W_t)}{\der t} 
     \le - m s^{2-\frac{2}{N}}L_{1,1}(W_t) \nonumber\\ 
     & \cdot \begin{cases}
     N \l((\a-\g_Z)^{2-\frac{2}{N}} - 2\g_Z \b^{2-\frac{2}{N}} \r) & t\le \tau \\
     \l(\a (\a-\g_Z)^{2-\frac{2}{N}} - 2\g_Z N \b^{2-\frac{2}{N}} \r)  & t \ge \tau . 
    \end{cases}
    \label{eq:overallLoss}
\end{align}
Suppose that inverse spectral gap $\g_Z$ is small enough such that the right-hand side of~\eqref{eq:overallLoss} is negative, see~\eqref{eq:simpler}. Using this  observation that $L_{1,1}(W_t)$ is decreasing in $t$ and by applying the Gronwall’s inequality to~\eqref{eq:overallLoss}, we arrive at Theorem~\ref{thm:MainSimple}.

\subsection{Derivation of (\ref{eq:evolveFast},\ref{eq:evolveSlow})}\label{sec:derivationEvolveFastSlow}

We begin with the detailed derivation of~\eqref{eq:evolveFast}. 
Let us repeat~\eqref{eq:evolveSimpleLemma2} for convenience:
\begin{align}
    & \frac{\der L_{1,1}(W_t)}{\der t} \nonumber \\
     & \le -  m N ((\a-\g_Z) s_Z)^{2-\frac{2}{N}} (s_t - u_t^\top Z_1 v_t)^2 
     \nonumber\\
      &  -2  \a m s_Z ((\a-\g_Z) s_Z)^{2-\frac{2}{N}}  (s_Z-u_t^\top Z_1 v_t)  \nonumber\\
      &  + m  N (\b s_Z)^{2-\frac{2}{N}} \g_Z |s_t - u_t^\top Z_1 v_t| (s_Z -u_t^\top Z_1 v_t) \nonumber\\
      &  + 2 m (\b s_Z)^{2-\frac{2}{N}}  s_{Z,2} (s_Z-u_t^\top Z_1 v_t).  \quad \text{(see \eqref{eq:evolveSimpleLemma2})} \nonumber
\end{align}
Recall also that~\eqref{eq:fastCvg} is in force. By ignoring the nonpositive term in the third line above, we arrive at
\begin{align}
         & \frac{\der L_{1,1}(W_t)}{\der t} \label{eq:evolveGeneralPrSimple}\\
     & \le  -  m N ((\a-\g_Z) s_Z)^{2-\frac{2}{N}} (s_t - u_t^\top Z_1 v_t)^2 \nonumber\\
      & \qquad +  m  N (\b s_Z)^{2-\frac{2}{N}} \g_{Z}  |s_t - u_t^\top Z_1 v_t| \nonumber\\
      & \qquad \qquad  \cdot \sqrt{s_Z(s_Z - u_t^\top Z_1 v_t)} \nonumber\\
     & \qquad + 2 m (\b s_Z)^{2-\frac{2}{N}}  \g_{Z}  s_Z (s_Z-u_t^\top Z_1 v_t).\nonumber
\end{align}
To obtain the first inequality in~\eqref{eq:evolveGeneralPrSimple},  note that $W_t\in \neigh_{\a,\b}(Z_1)$ by~\eqref{eq:WstaysNew} and, in particular, 
\begin{align}
u_t^\top Z_1 v_t \ge 0, \qquad \forall t\ge 0, \label{eq:innProdPos}
\end{align}
by definition of $\neigh_{\a,\b}(Z_1)$ in~\eqref{eq:neighWTwo}.
In turn,~\eqref{eq:innProdPos} implies that $s_Z - u_t^\top Z_1 v_t \le s_Z$.

We continue to bound the right-hand side of~\eqref{eq:evolveGeneralPrSimple} as
\begin{align*}
         & \frac{\der L_{1,1}(W_t)}{\der t}  \\
      & \le  -  m N ((\a-\g_Z) s_Z)^{2-\frac{2}{N}} (s_t - u_t^\top Z_1 v_t)^2 \nonumber\\
      & \qquad +  \frac{m  N}{\sqrt{2}}  (\b s_Z)^{2-\frac{2}{N}} \g_{Z}  (s_t - u_t^\top Z_1 v_t)^2    \nonumber\\
      & \qquad +  m (\b s_Z)^{2-\frac{2}{N}}  \g_{Z}  (s_t - u_t^\top Z_1 v_t)^2 = -m s_Z^{2-\frac{2}{N}} \cdot   \nonumber\\
      & \l(  (\a-\g_Z)^{2-\frac{2}{N}}  N  - \g_Z\b^{2-\frac{2}{N}} \l(\frac{N}{\sqrt{2}}+1\r) \r) \nonumber\\
      & \qquad \cdot (s_t - u_t^\top Z_1 v_t)^2 \le -m s_Z^{2-\frac{2}{N}}\cdot \nonumber\\
      & 
      \l(  (\a-\g_Z)^{2-\frac{2}{N}}  N  - \g_Z\b^{2-\frac{2}{N}} \l(\frac{N}{\sqrt{2}}+1\r) \r)\cdot  L_{1,1}(W_t).\nonumber
\end{align*}
The first inequality above uses~\eqref{eq:fastCvg} and the last inequality above uses~\eqref{eq:lossFastCvg}.

We next derive~\eqref{eq:evolveSlow}. Let us repeat~\eqref{eq:evolveSimpleLemma2} for convenience: 
\begin{align*}
     & \frac{\der L_{1,1}(W_t)}{\der t} \nonumber\\
          & \le -  m N ((\a-\g_Z) s_Z)^{2-\frac{2}{N}} (s_t - u_t^\top Z_1 v_t)^2 
     \nonumber\\
      &  -2\a    m s_Z ((\a-\g_Z) s_Z)^{2-\frac{2}{N}}  (s_Z-u_t^\top Z_1 v_t)  \nonumber\\
      &  + m  N (\b s_Z)^{2-\frac{2}{N}} \g_Z |s_t - u_t^\top Z_1 v_t| (s_Z -u_t^\top Z_1 v_t) \nonumber\\
      &  + 2 m (\b s_Z)^{2-\frac{2}{N}}  s_{Z,2} (s_Z-u_t^\top Z_1 v_t). 
\end{align*}      
Recall that now~\eqref{eq:slowCvg} is in force. By ignoring the  nonpositive term in the second line above, we then simplify the above bound as
\begin{align}
& \frac{\der L_{1,1}(W_t)}{\der t} \nonumber\\
     & \le 
        -  \a m  ((\a-\g_Z) s_Z)^{2-\frac{2}{N}} \cdot 2s_Z (s_Z-u_t^\top Z_1 v_t)  \nonumber\\
      &  + \frac{m  N}{\sqrt{2}} (\b s_Z)^{2-\frac{2}{N}} \g_Z \cdot 2 s_Z  (s_Z -u_t^\top Z_1 v_t) \nonumber\\
      &  +  m (\b s_Z)^{2-\frac{2}{N}} \g_Z \cdot  2 s_{Z} (s_Z-u_t^\top Z_1 v_t) \nonumber\\
      & = - m  s_Z^{2-\frac{2}{N}} \l( \a (\a-\g_Z)^{2-\frac{2}{N}} - \b^{2-\frac{2}{N}} \g_Z \l(\frac{N}{\sqrt{2}} + 1\r) \r) \nonumber\\
      & \qquad \cdot  2 s_{Z} (s_Z-u_t^\top Z_1 v_t) \nonumber\\
      & \le - m  s_Z^{2-\frac{2}{N}} \l( \a (\a-\g_Z)^{2-\frac{2}{N}} - \b^{2-\frac{2}{N}} \g_Z \l(\frac{N}{\sqrt{2}} + 1\r) \r) 
      \nonumber\\
      & \qquad \cdot L_{1,1}(W_t).
      \label{eq:slowRegimeEvolveSimplify}
\end{align}
To obtain the first inequality above, note that 
\begin{align*}
    \frac{1}{2} (s_t - u_t^\top Z_1 v_t)^2 & \le s_Z (s_Z - u_t^\top Z_1 v_t)
    \qquad \text{(see \eqref{eq:slowCvg})} \nonumber\\
    & \le s_Z^2 \qquad \text{(see \eqref{eq:innProdPos})},
\end{align*}
which, after rearranging, reads as 
\begin{align}
    |s_t - u_t^\top Z_1 v_t| \le \sqrt{2}s_Z.
\end{align}
The last inequality in~\eqref{eq:slowRegimeEvolveSimplify}  uses~\eqref{eq:slowCvgLoss}.

\subsection{Derivation of~(\ref{eq:conservative})}\label{sec:derConservative}
In the slow convergence regime in~\eqref{eq:slowCvg}, it holds  that 
\begin{align}
    \frac{1}{2}(s_t-u_t^\top Z_1 v_t)^2 & \le s_Z (s_Z - u_t^\top Z_1 v_t) 
    \qquad \text{(see \eqref{eq:slowCvg})}
    \nonumber\\
    & \le s_Z^2. \qquad \text{(see \eqref{eq:innProdPos})}
    \label{eq:conservativePr1}
\end{align}
On the other hand, we can also lower bound the first term  in~\eqref{eq:conservativePr1} as 
\begin{align}
    s_Z^2 
    & \ge \frac{1}{2}(s_t - u_t^\top Z_1 v_t)^2 
        \qquad \text{(see \eqref{eq:conservativePr1})}
    \nonumber\\
    & \ge \frac{s_t^2}{4} - \frac{1}{2}(u_t^\top Z_1 v_t)^2 \nonumber\\
    & \ge \frac{s_t^2}{4} - \frac{s_Z^2}{2},
        \label{eq:conservativePr2}
\end{align}
where the penultimate line above uses the inequality $(a-b)^2 \ge \frac{a^2}{2}-b^2$ for scalars $a,b$.  The last inequality above holds because $u_t,v_t$ are both unit-norm vectors and $s_Z$ is the only nonzero singular value  of $Z_1$, see~(\ref{eq:thinSVDrecall},\ref{eq:decomZ}), and thus $|u_t^\top Z_1 v_t| \le s_Z$.  

By rearranging~\eqref{eq:conservativePr2}, we find that 
$
    s_t > \sqrt{6} s_Z
$
implies the fast convergence regime  in~\eqref{eq:fastCvg}.

\end{document}